\documentclass[journal]{new-aiaa}
\usepackage[utf8]{inputenc}

\usepackage{graphicx}
\usepackage{amsmath}
\usepackage[version=4]{mhchem}
\usepackage{siunitx}
\usepackage{longtable,tabularx}
\title{Optimal Satellite Constellation Spare Strategy Using Multi-Echelon Inventory Control}

\author{Pauline Jakob \footnote{Formerly Graduate Research Assistant, Department of Aerospace Engineering, 104 S Wright St, Urbana, IL, 61801, USA.}}
\affil{University of Illinois at Urbana-Champaign, Urbana, IL, 61801, USA}
\author{Seiichi Shimizu\footnote{Advanced Technology R\&D Center.} and Shoji Yoshikawa \footnote{Advanced Technology R\&D Center.}}
\affil{Mitsubishi Electric Corporation, Amagasaki, 661-8861, Japan}
\author{Koki Ho\footnote{Assistant Professor, Department of Aerospace Engineering, 104 S Wright St, Urbana, IL, 61801, USA, AIAA Member.}}
\affil{University of Illinois at Urbana-Champaign, Urbana, IL, 61801, USA}

\usepackage{amssymb,amsmath}
\usepackage{lscape}
\usepackage[hyphens]{url}
\usepackage{hyperref}
\usepackage{setspace}
\usepackage{wasysym}
\usepackage{multirow}
\usepackage{array}
\usepackage{xcolor}
\newcolumntype{P}[1]{>{\centering\arraybackslash}p{#1}}
\begin{document}

\maketitle
\begin{abstract}
The recent growing trend to develop large-scale satellite constellations (i.e., mega-constellation) with low-cost small satellites has brought the need for an efficient and scalable maintenance strategy decision plan. Traditional spare strategies for satellite constellations cannot handle these mega-constellations due to their limited scalability in the number of satellites and/or frequency of failures. In this paper, we propose a novel spare strategy using an inventory management approach. We consider a set of parking orbits at a lower altitude than the constellation orbits for spare storage, and model the satellite constellation spare strategy problem using a multi-echelon (s,Q)-type inventory policy, viewing the Earth’s ground as a supplier, the parking orbit spare stocks as warehouses, and the in-plane spare stocks as retailers. The accuracy of the proposed analytical model is assessed using simulations via Latin Hypercube Sampling. Furthermore, based on the proposed model, an optimization formulation is introduced to identify the optimal spare strategy, comprising the parking orbits' characteristics and all locations' policies, to minimize the maintenance cost of the system given performance requirements. The proposed model and optimization method are applied to a real-world satellite mega-constellation case to demonstrate their value.

\end{abstract}

\section*{Nomenclature}

{\renewcommand\arraystretch{1.0}
\noindent\begin{longtable*}{@{}l @{\quad=\quad} l@{}}
$cap_{\text{launch}}$ & Launch capacity (number of possible satellites per rocket), in units of satellites\\
$D_{\text{plane}}$ & Demand for in-plane spares, in units of satellites\\
$D_{\text{parking}}$ & Demand for parking spares, in units of batches $Q_{\text{plane}}$\\
$ES_{\text{plane}}$ & Expected number of backorders for in-plane spares over a replenishment cycle, in units of satellites\\
$ES_{\text{parking}}$ & Expected number of backorders for parking spares over a replenishment cycle, in units of batches $Q_{\text{plane}}$\\
$f_{\text{plane}}$ & Probability density function of the lead time to the constellation orbits, in units of days$^{-1}$\\
$f_{\text{parking}}$ & Probability density function of the lead time to the parking orbits, in units of days$^{-1}$\\
$h_{\text{plane}}$ & Altitude of the constellation orbits, in units of kilometers\\
$h_{\text{parking}}$ & Altitude of the parking orbits, in units of kilometers\\
$i$ & Inclination of the constellation orbits, in units of degrees\\
$k_{\text{Q,parking}}$ & Order batch size for parking spares, in units of batches $Q_{\text{plane}}$\\
$k_{\text{s,parking}}$ & Safety stock for parking spares, in units of batches $Q_{\text{plane}}$\\
$\lambda_{\text{sat}}$ & Failure rate of a satellite, in units of failures per year\\
$\lambda_{\text{plane}}$ & Demand rate for in-plane spares, in units of satellites per day\\
$\lambda_{\text{parking}}$ & Demand rate for parking spares, in units of batches $Q_{\text{plane}}$ per day\\
$m_{\text{dry}}$ & Dry mass of the satellites, in units of kilograms\\
$m_{\text{fuel}}$ & Fuel mass required for a Hohmann transfer (from a parking orbit to a constellation orbit), in units of kilograms\\
$\mu_{\text{launch}}$ & Mean time between launches to the parking orbits, in units of days \\
$N_{\text{days}}$ & Number of days per year\\
$N_{\text{fail,plane}}(\tau)$ & Demand for in-plane spares during a lead time $\tau$,  in units of satellites\\
$N_{\text{fail,parking}}(\tau)$ & Demand for parking spares during a lead time $\tau$, in units of batches $Q_{\text{plane}}$\\
$N_{\text{plane}}$ & Number of constellation orbital planes\\
$N_{\text{parking}}$ & Number of parking orbital planes\\
$N_{\text{sats}}$ & Number of operational satellites per orbital plane in the constellation\\
$P_{\text{av}}$ & Parking orbit availability\\
$P(i^{th})$ & Probability of getting supply from the $i^{th}$ closest parking orbit\\
$p_{\text{launch}}$ & Launch cost, in units of million US\$ per launch \\
$p_{\text{sat}}$ & Manufacturing cost of unit satellite, in units of million US\$ per satellite\\
$p_{\text{holding}}$ & Annual holding cost for each satellite, in units of million US\$ per satellite per year\\
$p_{\text{launch,full}}$ & Cost of a full rocket launch (for $cap_{\text{launch}}$), in units of million US\$ per launch\\
$p_{\text{launch,unit}}$ & Cost of a unique satellite rocket launch (for one satellite only), in units of million US\$ per launch\\
$pt_{\text{launch}}$ & Order processing time for launch, in units of days\\
$Q_{\text{plane}}$ & Batch size for in-plane spares, in units of satellites\\
$Q_{\text{parking}}$ & Batch size for parking spares, in units of satellites\\
$\rho_{\text{plane}}$ & Order fill rate for in-plane spares\\
$\rho_{\text{parking}}$ & Order fill rate for parking spares\\
$s_{\text{plane}}$ & Reorder point for in-plane spares, in units of satellites\\
$s_{\text{parking}}$ & Reorder point for parking spares, in units of satellites\\
$\overline{Stock_{\text{plane}}}$ &Mean stock level of in-plane spares, in units of satellites\\
$\overline{Stock_{\text{parking}}}$ & Mean stock level of parking spares, in units of batches $Q_{\text{plane}}$\\
$T_{\text{plane}}$ & Lead time to the constellation orbits, in units of days\\
$T_{\text{parking}}$ & Lead time to the parking orbits, in units of days\\
$v_\text{exhaust}$ & Effective exhaust velocity of the satellite thrusters, in units of kilometers per second
\end{longtable*}}

\section{Introduction}

\begin{figure}[ht!]
\centering
\includegraphics[width = 12cm]{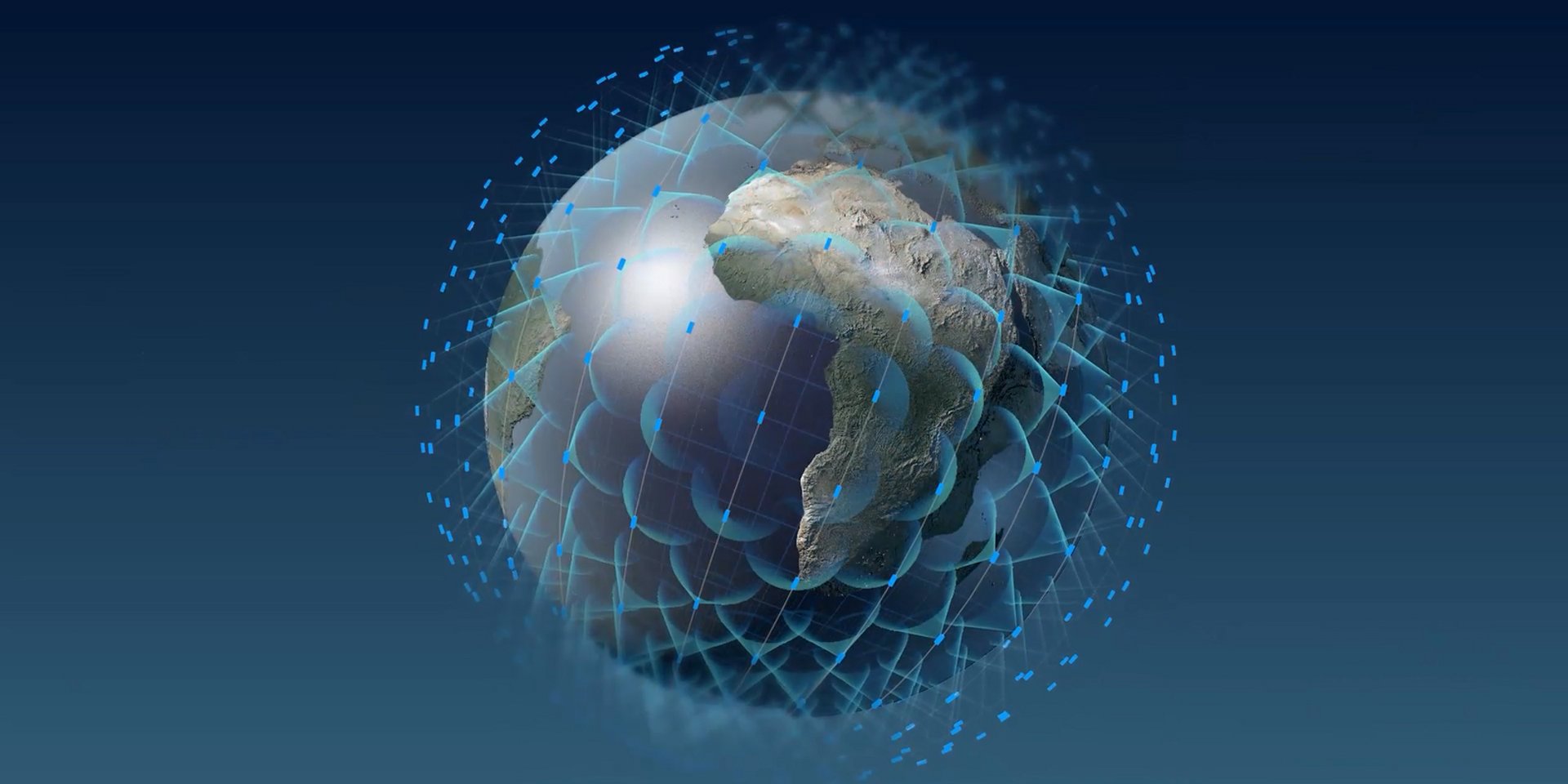}
\caption{OneWeb satellite constellation \cite{OneWebPic}}
\label{OneWeb}
\end{figure}

The trend for satellite constellations has been growing since their first establishment about twenty years ago, in May 1997 for Iridium and February 1998 for Globalstar. Various studies have been performed with their focuses on optimization of satellite constellation design \cite{ferringer2006satellite}\cite{ferringer2007efficient}\cite{budianto2004design}\cite{bandyopadhyay2016review}. More recently, new constellations comprising a tremendous number of small satellites (i.e., mega-constellation) have been considered to respond to the explosive demand for telecommunication services. For example, OneWeb is setting up a 900-satellite constellation in a low-Earth orbit (LEO) to provide broadband services (see Fig. \ref{OneWeb}) \cite{OneWebPic}, while SpaceX is planning a mega-constellation of nearly 12,000 interlinked broadband Internet satellites  \cite{spacex}\cite{fccmega}. 

In order to ensure the prosperity of the providers, guaranteeing a high level of customer satisfaction is paramount. Indeed, as discussed by Diekelman \cite{diekelman1998design}, the satellite failure mitigation can take a few days to several weeks, and the impact of a failed satellite can affect not only the current lost revenue but also the reputation of the system, and thus its future revenue. Therefore, it becomes vital to maintain the operational state of the system and secure a minimum availability to provide the offered services by avoiding outages. In the case of Iridium, for instance, twenty of the original satellites launched have required replacement, and spare satellites represent a substantial part of the constellation with about 30\% of the original fleet \cite{iridium}. As the trend for satellite mega-constellations grows, a large number of satellite failures can be expected from future mega-constellations, and a steady replacement strategy has to be established to maintain the service level.

Existing satellite constellation spare strategies are not effective for large-scale small satellite constellations. Traditional spare strategies include having some ground spare satellites to replace the failed satellites using on-demand launch, or having a few active or inactive spare satellites in every orbital plane for redundancy \cite{palmerini1998hybrid}.
Although these approaches were acceptable for small-scale constellations with large and highly reliable satellites (i.e., infrequent failures), they are not effective for mega-scale small satellite constellations, where each satellite tends to display less redundancy and thus less reliability to favor cost efficiency. Indeed, using only in-plane spares could result in needs for a large number of spare satellite units per orbital plane, thus involving a very high spare strategy cost. On the contrary, launching spare satellites on demand is a risky strategy, given the uncertainties in launch time schedule and the high satellite failure rate. Moreover, the launch of spare satellites itself can be problematic as typical rockets load tens of small satellites (e.g., 150 kg per satellite for OneWeb's constellation \cite{onewebBuilt}) per launch leveraging the batch launch discount (i.e., the cost-saving effects by launching many satellites in one rocket); we cannot provide on-demand launches for every spare satellite at a low cost. 
Some companies have foreseen the replenishment of their constellation-to-be, such as OneWeb, who signed a contract with Virgin Galactic to use their LauncherOne vehicle to haul up one satellite at a time. Yet, the solution adopted by OneWeb would raise the spare launch cost to be approximately seven times higher than that of a nominal satellite launch \footnote{The company targets about US\$10 million per resupply mission using the LauncherOne vehicle, whereas its contract for the initial constellation deployment with Arianespace values the launch of 700 satellites for US\$1 billion, and thus US\$1.43 million per satellite, packing 32 to 36 satellites per launch \cite{onewebLaunch}. }. Therefore, it is still beneficial if we could optimally take advantage of the batch launch discount, which was not possible in the traditional approaches.
There is a growing demand to have an automatic and scalable decision making and planning strategy under the uncertainty of satellite failures, in order to ensure the maintenance of the system. \footnote{Maintenance in this paper refers to the replacement of failed satellites to maintain the constellation as a system.} 

This paper offers a novel and unique design technique that is scalable for mega-scale satellite constellation replacement strategies leveraging inventory management methods. 
Our solution is to incorporate a set of parking orbits at a lower altitude than the constellation orbits to save on launch cost, and optimize the spare strategy as a supply chain between the Earth's ground (supplier), the parking orbits (warehouses), and the constellation orbits (retailers). 
A multi-echelon inventory control system is considered, under stochastic demands and lead times, comprising one supplier (Earth's ground), multiple warehouses (the parking orbits with parking spares), and multiple retailers (the constellation orbits with in-plane spares). An $(s,Q)$-type inventory policy is considered so that the system can optimally leverage the batch launch discount. An analytical model for the constellation spare strategy is developed in this paper, and an optimization formulation is introduced to optimize the spare strategy, minimizing the maintenance cost of the constellation. 
The developed optimization formulation can quickly approximate the optimal spare strategy without relying on computationally costly simulations; if necessary, the resulting optimized strategy can be further analyzed with high-fidelity simulations.

Although this paper mainly focuses on satellite systems, the general model developed in this paper also extends the existing inventory management literature. 
The interesting property of our problem is the specific interactions between the different levels of inventory on demand, lead times, and supply allocation.
Particularly, our problem is unique in that its multiple warehouses (the parking orbits) drift over time with respect to multiple retailers (the constellation orbits) due to orbital mechanics effects, and the retailers choose the closest (defined as the minimum waiting time in our context) available warehouse at the time of delivery. 
The general framework allowing retailers to get supplies from different warehouses can provide flexibility to avoid, or at least reduce, stock-out times. 
The accuracy of the analytical model developed in this paper is assessed using simulations, and the model is then leveraged for optimization of the spare strategy.
The proposed model and optimization formulation are applied to a real-world satellite mega-constellation case to demonstrate their value.

The remainder of the paper is organized in the following way. Section \ref{litreviewsection} presents an overview of the related literature from both the optimal satellite constellation spare strategy and supply chain model points of view. Section \ref{preliminariessection} provides the reader with preliminaries about the general theory of satellite constellations and inventory management useful for the understanding of the model further developed in Section \ref{modelsection}. Section \ref{validationsection} assesses the accuracy of the developed analytical model using simulations. In Section \ref{optimizationsection}, the optimization of the spare strategy is presented, and Section \ref{examplessection} provides a case study for the maintenance of a LEO communication satellite mega-constellation along with a sensitivity analysis for different satellite failure rates. Finally, Section \ref{conclusionsection} concludes the paper.

\section{Literature Review}\label{litreviewsection}
The literature regarding modeling of satellite constellations and their spare strategies is very sparse. Different solutions have been examined to ensure the replacement of failed satellites in orbit for such constellations. Lang and Adams \cite{lang1998comparison}, Lansard and Palmade \cite{lansard1998satellite}, Palmade et al. \cite{palmade1998skybridge} and Cornara et al. \cite{cornara1999satellite} all proposed global constellation design including analysis of their replacement strategies, choosing between distinct spare strategies including ground spares, parking orbit spares, in-plane spares, and overpopulated planes. 
However, no mixtures of each strategy have been considered, leaving the decision makers little flexibility in spare strategies.
Also, the complexity of such systems often leads authors to use simulations to evaluate the satellite reliability or constellation availability over time. However, the use of Monte Carlo simulations \cite{cornara1999satellite} can result in computational inefficiency, especially in the case of mega-scale constellations.

Other proposed models handled the simulation issue by adopting an analytical point of view and represented the satellite constellations by an exhaustive number of states; however, most of these models have significant scalability issues. Ereau and Saleman \cite{ereau1996modeling} approached the availability issue of satellite constellations using Petri nets, but in order to properly incorporate the use of time, the analytical results would still face the issue of state space explosion. Sumter \cite{sumter2003optimal} established an analytical model to find an optimal satellite replacement policy by the means of finite-horizon Markov decision processes, minimizing the expected monetary and opportunity costs of maintaining the constellation. The author limits the state explosion issue raised by Ereau and Saleman by setting several assumptions regarding satellites and their operation, such as zero launch lead time and only considering ground spares. Those suppositions can be questionable and Sumter also recognizes the limitations in the work. Furthermore, the number of states considered for the solution regarding the size of the constellation still remains very large especially for a mega-constellation, with, for example, 4,608 states for a constellation comprising nine satellites. Kelley and Dessouky \cite{kelley2004minimizing} also used Markov models to evaluate the life cycle cost of a satellite system comprising acquisition, replenishment, and operations costs, linked to a performance model to assess the availability of the service. Again, this type of modeling leads to state explosion as the size of the constellation increases, and thus is not scalable to planned mega-constellations.

There have been very few attempts to model the orbiting satellite constellation spare strategy problem using an inventory management approach. Dishon and Weiss \cite{dishon1966communications} originally analyzed the problem of satellite replenishment from a simple satellite level perspective and solved it using a classical $(N,M)$ inventory system. Their solution would consider the total number of functional satellites in a given system, and when the latter falls from $M$ to $N$, replenishment launches are initiated to repopulate the system up to level $M$. An optimal policy was derived using a number-of-satellites-launched-over-time cost function. However, the considered inventory model was very simple and presented several limitations: the replenishment up to a level $M$ does not allow consistent launch planning over time; it cannot reflect the reality of batch launch discount and it does not explicitly consider the use of parking orbits. These limitations make the proposed strategy ineffective for large-scale satellite constellations.

Although very few authors developed a satellite constellation replenishment policy leveraging inventory management techniques, the general problems of spare parts inventory control and supply chain management have been studied widely in the literature. Many mathematical models have been proposed over time for supply chain inventories. Multi-echelon systems are particularly interesting for the purpose of satellite constellation spare strategies and different papers have tried to grasp the complex interactions between different levels of such systems subject to various features. A detailed review can be found in Ref. \cite{gumus2007multi}. In this impetus, Ganeshan \cite{ganeshan1999managing} followed the work of Deuermeyer and Schwarz \cite{schwarz1981} and developed a model for multi-level inventory comprising multiple retailers, one warehouse and multiple identical suppliers while taking advantage of order splitting policies. Various applications of multi-level inventory policies can also be found in the literature. Costantino et al. \cite{costantino2013multi} presented an example of spare parts allocation using multi-echelon inventory control applied to the aeronautical industry, a very demanding sector in terms of availability requirements, while Caglar et al.\cite{caglar2004two} developed a continuous review, base stock policy for a two-echelon, multi-item spare parts inventory system for electronic machines. However, no model has been proposed and studied to address our unique challenge in the satellite constellation spare strategy, which requires multiple warehouses drifting over time, all able to resupply all the retailers and with stochastic demands at the retailers. 

In order to address this significant literature gap, our approach regarding the spare strategy for satellite constellations aims at concurrently considering different levels of spare satellites in the system, including ground spares, parking spares, and in-plane spares, and optimizes the whole supply chain using an analytical model with no need for simulations.

\section{Preliminaries}\label{preliminariessection}
The analysis in this paper builds upon concepts from two different fields: satellite constellations and inventory management. This section provides the readers with an appropriate description of the enabling notions needed to understand the concepts of satellite constellations and inventory management in the context of this paper.

\subsection{Satellite Constellations}
This subsection presents the theory of orbit perturbations, orbital transfers, and satellite constellations. Only the key elements relevant to this paper are explained here, and further theory can be found in Ref. \cite{prussing1993orbital}.

\subsubsection{Orbit perturbations} \label{Orbitpert}

To analyze the satellite constellation, we need to model the orbit of satellites orbiting the Earth. Only circular orbits are explored in this paper. The classical two-body orbital dynamics relies on the approximation that the Earth is a point mass; however, various factors can cause perturbations to the motion of satellites in reality. Two largest perturbations affecting a satellite's motion about the Earth are the atmospheric drag and the effects of Earth's oblateness. At the altitudes considered in this research ($\geq 700 \mathrm{km}$), atmospheric drag is considered to have negligible effects on the motion; however, the effects of Earth's oblateness are not negligible.

The oblateness of the Earth causes the irregularity in the gravitational field: the mass spinning creates an extra bulge around the equator, further causing perturbations to the satellite's orbital motion. This oblateness is characterized by a constant, $J_2 = 0.00108263$, contributing to a perturbing acceleration and disturbing the orbital elements. One of the principal effects of the Earth oblateness disturbance that is relevant to our research is to cause the right ascension of the ascending node (RAAN) $\Omega$ of an orbit to drift over time, with a rate depending only on the semimajor axis $a$ of this particular orbit and its inclination $i$:

\begin{equation}
\frac{d\Omega}{dt} = - \frac{3 n R_{\text{Earth}}^2J_2}{2a^2}\cos{i}
\label{drift}
\end{equation}

\noindent where $n = \sqrt{\frac{\mu}{a^3}}$ is the mean motion of the satellite, $\mu$ is the standard gravitational parameter of the Earth, and $R_\text{Earth}$ is the (mean) radius of the planet Earth. Note that $a$ is a function of the altitude; therefore this change of RAAN depends on the altitude of the orbit.

\subsubsection{Orbital transfer} \label{Orbittransfer}
In order to deliver a spare from one orbit (e.g., a parking orbit) to another orbit (e.g., a constellation orbit), we need to consider an orbital transfer. In this study, we consider Hohmann transfers, a common fuel-efficient type of chemical transfer for co-planar circular orbits. Out-of-plane maneuvers are excluded in this paper due to their cost inefficiency. In a Hohmann transfer, the cost of the transfer is evaluated through the mass of fuel $m_{\text{fuel}}$ required to perform the transfer, which depends on the velocity variation $\Delta V_\text{Hohmann}$ needed to move the satellite from an orbit altitude to another, the effective exhaust velocity $v_\text{exhaust}$ of the thruster, and the dry mass $m_{\text{dry}}$:

\begin{equation}
\label{mfuel}
m_{\text{fuel}} = m_{\text{dry}}(e^{\frac{\Delta V_\text{Hohmann}}{v_\text{exhaust}}}-1)
\end{equation}
where $\Delta V_\text{Hohmann}$ can be calculated based on the radius (i.e., semimajor axis) of the initial orbit $a_0$ and the final orbit $a_1$ ($a_1>a_0$ in the context of this paper) as follows \footnote{\label{footnoteprussing} See \cite{prussing1993orbital} for more details about Hohmann transfer calculations.}:

\begin{equation}
\label{deltaV}
\Delta V_\text{Hohmann} = \sqrt{\frac{\mu}{a_0}}(\sqrt{\frac{2a_1}{a_0+a_1}}-1)+ \sqrt{\frac{\mu}{a_1}}(1-\sqrt{\frac{2a_0}{a_0+a_1}})
\end{equation}

The time of flight of such a Hohmann transfer corresponds to half a period of the transfer ellipse of the semimajor axis $\frac{a_0+a_1}{2}$:
\begin{equation}
\label{tof}
TOF_H = \pi \sqrt{\frac{(a_0+a_1)^3}{8\mu}}
\end{equation}

\subsubsection{Satellite constellation}
A satellite constellation is a set of satellites working together in order to provide a service. When the number of satellites comprised in the system becomes very large, we denote it as a \textit{mega-constellation}. The well-known Walker Delta pattern constellation  \cite{1984JBIS...37..559W} is considered in this paper. In this configuration, the total number of satellites is allocated to $N_{\text{plane}}$ circular orbital planes (i.e., referred to as the constellation orbits), such that there are $N_{\text{sats}}$ satellites per plane. All constellation orbits share the same altitude $h_{\text{plane}}$ and the same inclination $i$, and their RAANs $\Omega$ are distributed such that the planes are equally spaced ($\Omega_{k\text{-th plane}} = (k-1)*\frac{2\pi}{N_{\text{plane}}}$). This strategy is of particular interest to preserve the geometry of the system, as all satellites would endure approximately the same orbit perturbations. In other words, all satellites in the constellation would experience the same RAAN drift rate.  

Therefore, two constellations with the same inclination but different altitudes (e.g., the constellation orbits and the parking orbits) would have two distinct nodal shift rates $\frac{d\Omega}{dt}$ and thus we can observe a relative RAAN drift between them. The spare strategy model utilized in this paper takes advantage of this specific relative RAAN drift between the constellation orbits and the parking orbits, the latter of which are located at a lower altitude.

\subsection{Inventory Management \label{secinventory}}
Inventory management considers the flow of products (e.g., spare parts in our context) in a supply chain and enables delivery of a better service. It encompasses the relations between all levels of inventory, from suppliers to warehouses and to retailers. Inventory control is of primary importance, especially for stochastic demands and lead times. In this subsection, the specific (s,Q)-policy is introduced along with its characteristic features such as the replenishment cycles, backorders, and mean stock level. Note that the model presented here assumes that stock-outs happen rarely and thus are negligible for calculation of the stock level; this assumption is common in the literature \cite{ganeshan1999managing} and is also reasonable for our application as discussed later.

\subsubsection{(s,Q)-policy}
\begin{figure}[ht!]
\centering
\includegraphics[width = 14cm]{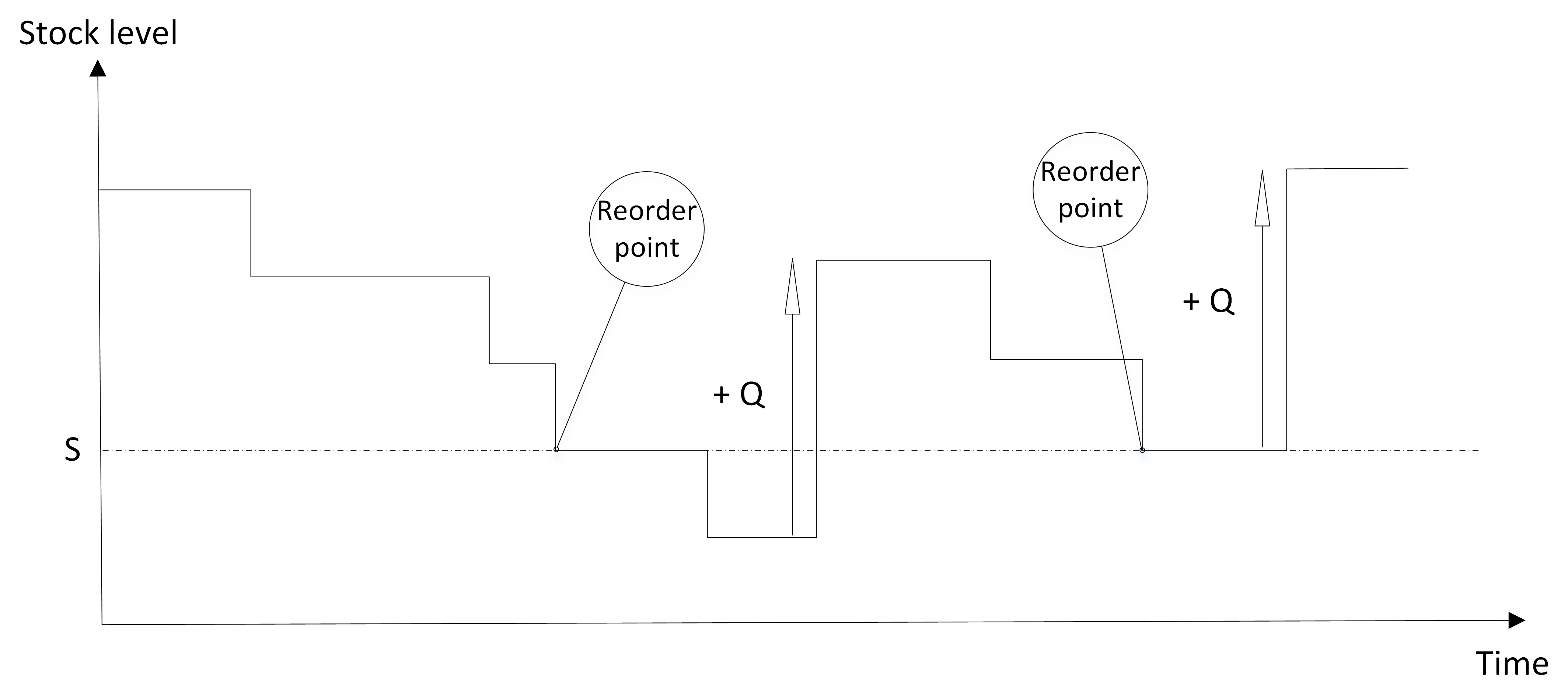}
\caption{Representation of replenishment cycles}
\label{cycle}
\end{figure}

All the facilities considered in the model are assumed to follow a continuous $(s,Q)$-type inventory policy. This particular policy is chosen because it enables optimization of the order quantity $Q$, unlike other policies such as $(R,S)$ or $(s,S)$ policies, so that we can maximize the batch launch discount. In the $(s,Q)$ inventory policy, each facility (e.g., warehouse, retailer) holds an inventory of the spare stock, and when a stock level drops to or below $s$ available units, an order of batch size $Q$ is placed to its attached supplier. The parameters $s$ and $Q$ can be optimized.
The model presented in this paper focuses on the study of replenishment cycles, each of which contains a replenishment of $Q$ units.  
Fig. \ref{cycle} illustrates replenishment cycles from a stock point-of-view. 

\subsubsection{Backorders} \label{EStau}
The model takes the situation of backorders into consideration in order to evaluate the efficiency of the policy. When a demand cannot be met by on-stock spare units, it is backordered. The next spares supply has to satisfy this backordered demand first upon arrival. It is important to be able to evaluate the short units (i.e., backorders) that the different facilities would be facing over replenishment cycles and have the means to control them. Knowing that the replenishment phase starts when the stock level drops to or below $s$, the expected backorders per cycle $ES$ for a lead time $\tau$ become the demand exceeding $s$ units during the time $\tau$, which can be derived from the probability distribution of failures in Eq.\ref{expectedshortage} \cite{ganeshan1999managing}\cite{costantino2013multi}.

\begin{equation}
\label{expectedshortage}
ES = \sum_{k \ge s+1} \ (k - s) \ P_{\tau}(D = k)
\end{equation}

\noindent where $P_{\tau}(D = k)$ is the probability of having $k$ units of demand during a lead time $\tau$. 

In order to manage the number of backorders that a facility would face, we introduce the notion of order fill rate $\rho$, which is the percentage of demand that is satisfied from the available stock during a cycle. The order fill rate can be found using Eq. \ref{fillrate}.

\begin{equation}
\label{fillrate}
\rho=1-\frac{ES}{Q}
\end{equation}
\noindent The order fill rate is linked to the performance of the replenishment policy at a facility. The optimal design is chosen so that the global multi-echelon spare system meets performance requirements.

\subsubsection{Mean stock level \label{subsecstock}}
\begin{figure}[ht!]
\centering
\includegraphics[width = 10cm]{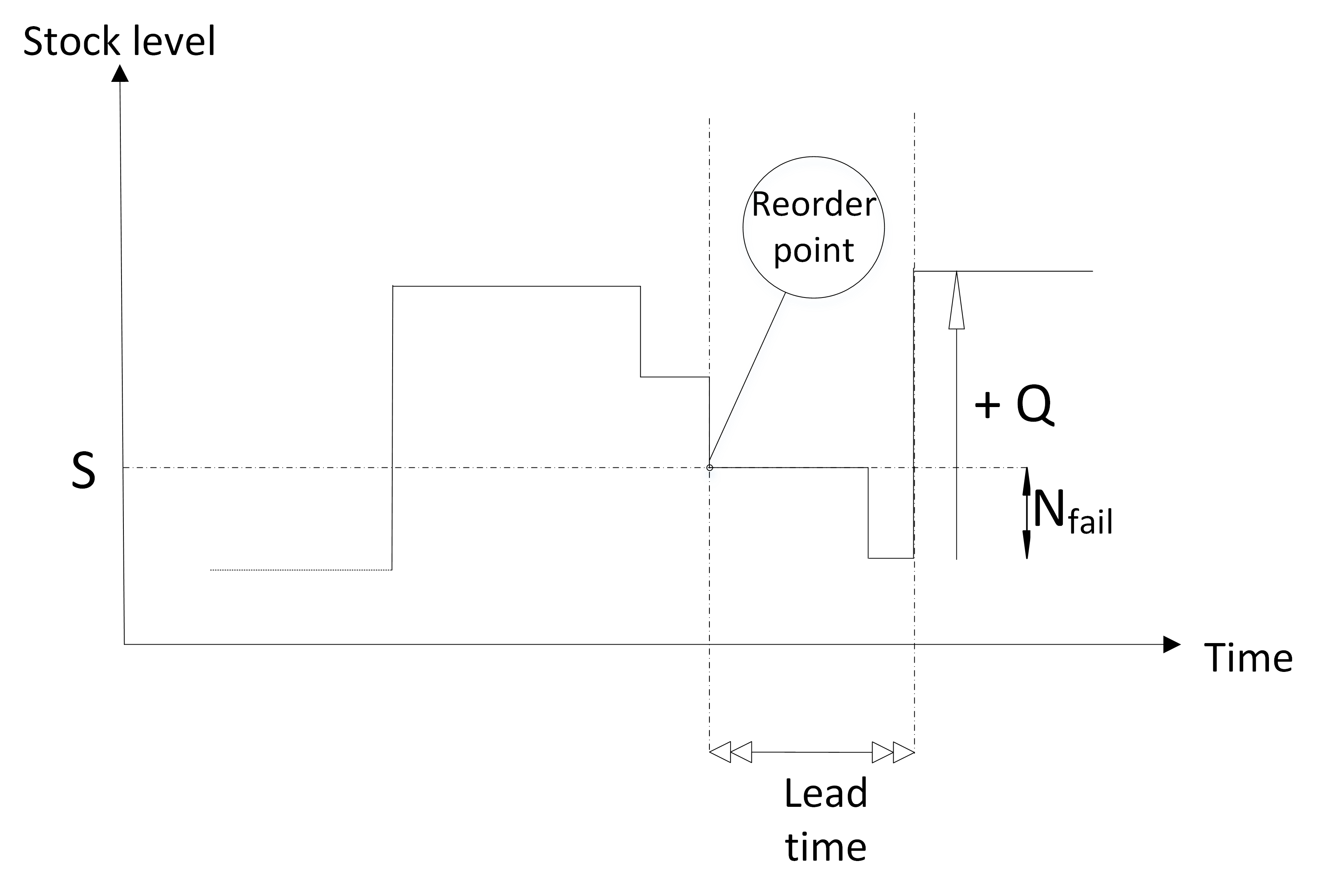}
\caption{Illustration of the stock level}
\label{cyclezoom}
\end{figure}

It is of particular interest to know the mean stock level at each facility to be able to further derive holding costs. From Fig. \ref{cyclezoom}, the stock level is comprised between $(Q+s-N_{\text{fail}}(\tau))$ and $(s-N_{\text{fail}}(\tau))$, where $N_{\text{fail}}(\tau)$ is the number of failures during a lead time $\tau$. Thus, assuming a linear continuous stock level drop, the mean stock level would be $(\frac{Q}{2} + s - N_{\text{fail}}(\tau))$.  Furthermore, the continuity correction factor of $1/2$ needs to be added to adjust the difference between the real discretized stock level drops and the assumed linear continuous stock level drops. \footnote{Note that this continuity correction factor can be negligible for a large stock level, but it can be important for our application with a relatively small stock level.} Eq. \ref{meanstock} gives the resulting average stock level \cite{jensen2003operations}:
\begin{equation}
\label{meanstock}
\overline{Stock}= \frac{Q}{2} + s - N_{\text{fail}}(\tau) +\frac{1}{2}
\end{equation}

\section{Model Formulation}\label{modelsection}

\subsection{Overview of the Model}

\begin{figure}[ht!]
\centering
\includegraphics[width=10cm]{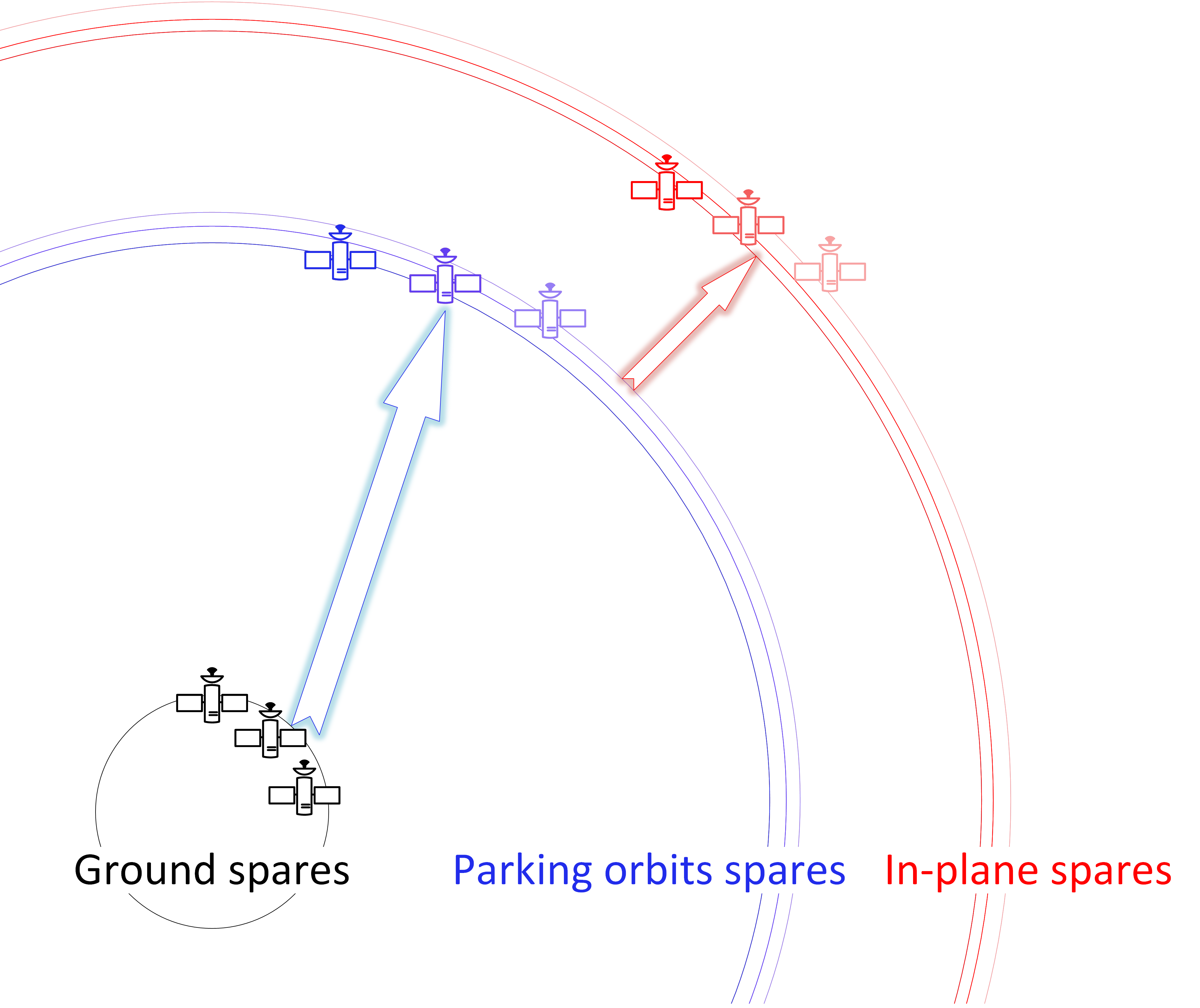}
\caption{Overview of the multi-level spare strategy for a satellite constellation}
\label{strategy}
\end{figure}

The aim of the model is to provide a replenishment strategy for the spare parts of a satellite constellation and establish a criterion to evaluate maintenance strategy performances. As presented by Cornara et al. \cite{cornara1999satellite}, different spare strategies exist to ensure the maintenance of the constellation (see Table \ref{strategies}). To provide more flexibility in the design of the spare strategy for satellite constellations, this paper introduces a mixed-strategy with multiple levels of spares, taking advantage of each approach. A visual representation of the strategy is given in Fig. \ref{strategy}.

\begin{table}[ht!]
\centering
\renewcommand{\arraystretch}{1.3}
  \caption{Different possible spare strategies and their approximate replacement time proposed by Ref. \citep{cornara1999satellite}}
\begin{tabular}{cc}
  \hline
  \hline
  Strategy & Replacement time\\
  \hline
  Overpopulation  & No time delay for replacement \\
  In-plane spares & 1-2 days\\
  Parking spares  & 1-2 months \\
  Ground  & A few months to 1 year \\
  \hline
  \hline
  \end{tabular}  
  \label{strategies}
\end{table}

The first level of spares is the constellation's in-plane spares. The paper does not distinguish between active (overpopulation strategy) or inactive (in-plane strategy) spares and lets this choice to the reader. When a satellite failure occurs in one of the constellation orbits, and if a spare part is in stock in that orbital plane, the failed satellite is replaced using available in-plane spares. This first level allows the constellation to avoid outages with little to no time delay to replace a failed satellite. 

The second level of spares is parking spares. It consists of spare satellites placed in a lower altitude orbit and at the same inclination as the constellation orbits, and are available to transfer to the in-plane spare stocks using orbital maneuvers. Note that there can be one or multiple parking orbits (i.e., multiple orbital planes); all parking orbits are circular, share the same altitude and inclination, and have their RAANs equally spaced. When the spare stock level of the in-plane locations reaches a critical level, it places an order to the parking orbits to be resupplied with spare satellites. Having this second level of spares available provides the orbital planes with the possibility to replenish their spare stocks within a relatively short amount of time (see Table \ref{strategies}), and can thus reduce the number of spares needed in each constellation orbit. In addition, since the parking orbits can replenish the spare stocks of any constellation orbit, they increase the flexibility of the supply chain.

Finally, the last level of spares is ground spares, i.e., spare satellites on the Earth's ground, which are considered to be always available to replenish the parking orbits thanks to the fast manufacturing assembly line that is achievable nowadays for satellite constellations \cite{onewebBuilt}. Whenever a parking orbit reaches its critical stock level, it places an order to the ground spare stock to schedule a rocket launch to replenish its stock. 

\begin{figure}[ht!]
\centering
\includegraphics[width=16cm]{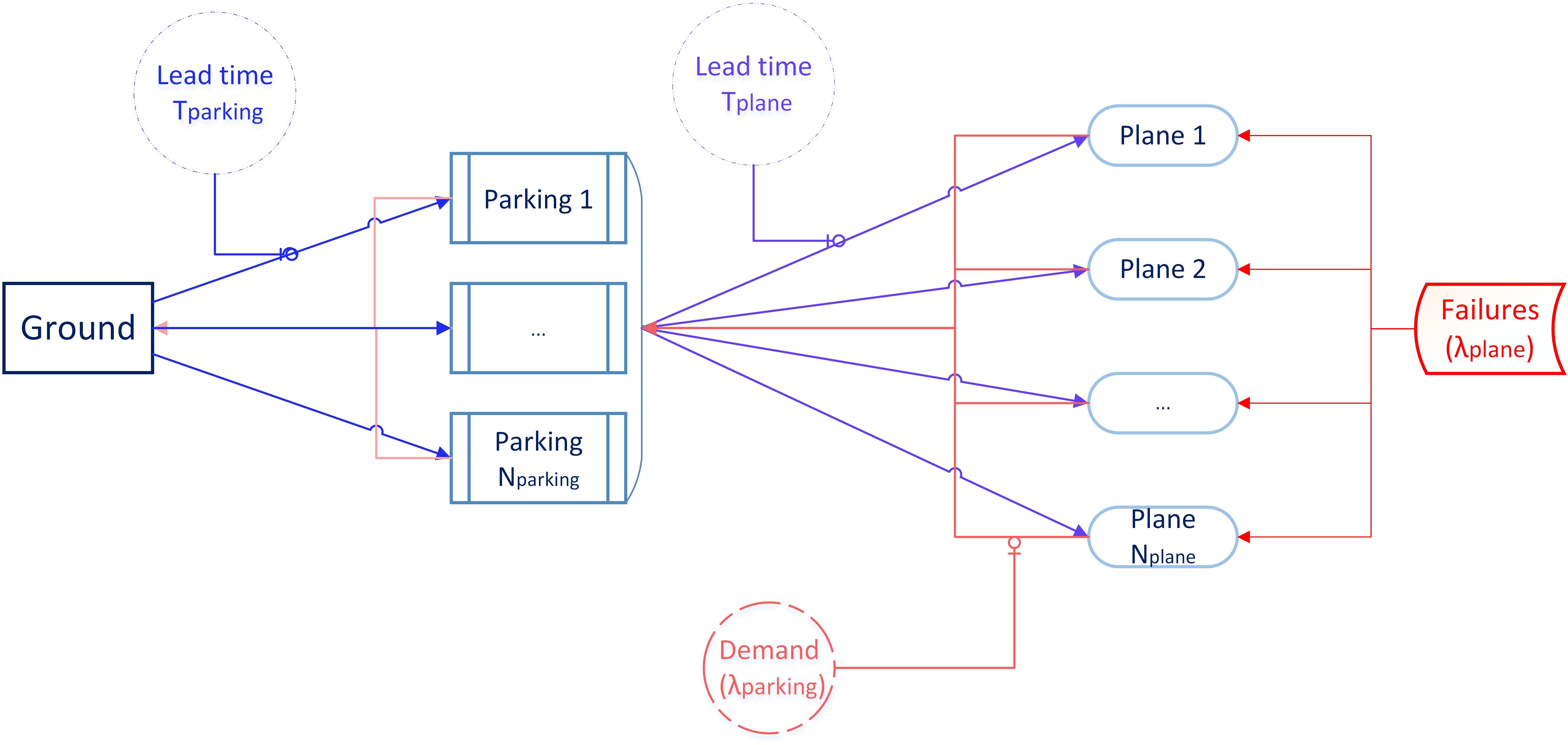}
\caption{Proposed multi-echelon inventory model for a satellite constellation}
\label{model}
\end{figure}

All levels of spare locations together are considered as a multi-echelon inventory system, with stochastic demands associated with satellite failures, and stochastic lead times for both types of replenishment (from the ground to the parking orbits, and from the parking orbits to the constellation orbits). Fig. \ref{model} captures the interactions between the different levels of inventory.

The remainder of this section is organized as follows. Section \ref{assumptions} introduces the general assumptions used in our model. Section \ref{inplanemodel} and Section \ref{parkingmodel} are symmetric, as they introduce the inventories of in-plane spares and parking spares, respectively. Finally, Section \ref{totcostmodel} describes the cost model used to evaluate the spare strategy, which will be used for the optimization in Section \ref{optimizationsection}.

\subsection{Model Assumptions} \label{assumptions}

The following presents a summary of the general assumptions of our model. Other assumptions are discussed as the model is introduced in more detail.

\begin{itemize}

\item Spare parts located in the first echelon (in-plane spares) are considered to be immediately available to replace a failed satellite unit. This postulate is true in the case of an overpopulated strategy; however, in case of spare satellites located in a slightly different plane to avoid collisions, the process of replacement can take up to 2 days. This delay is not considered in the model.

\item The constellation's in-plane spare stocks get supplies from the closest (i.e., minimum waiting time) available parking orbit's spare stock. In order to allow flexibility in the spare replacement flow, we allow any parking orbit to potentially resupply any orbital plane's stocks. When a constellation orbit's in-plane stock reaches the re-order point ($s$-level), an order is placed to all parking orbits jointly and the spares batch is supplied from the closest parking orbit with spare availability at the time of the order. 

\item Supply from the ground can be delivered only to a unique parking orbit. Indeed, as stated by Lang and Adams \cite{lang1998comparison}, using a single rocket launch to supply different orbital planes can turn out to be very inefficient. 

\item To facilitate the tracking of the orders, an order is allowed to be processed only when no previous order is already in transit. 

\item Stock-outs are assumed to happen rarely and thus are negligible for calculation of the stock levels. This assumption is reasonable for our optimal spare strategies.

\item  As the spares have to be transferred by batches both from the Earth's ground to the parking orbits and from the parking orbits to the constellation orbits, the order quantity and re-order point at the parking orbits are assumed to be multiples of the batch size $Q_{\text{plane}}$ of in-plane spares:

\begin{equation}\label{batcheq}
\left\{
\begin{array}{rl}
  Q_{\text{parking}} &= k_{\text{Q,parking}} \ Q_{\text{plane}} \\
s_{\text{parking}}&= k_{\text{s,parking}} \ Q_{\text{plane}}
\end{array}
\right.
\end{equation}

\end{itemize}

\subsection{In-plane Spares Inventory Model}\label{inplanemodel}
This subsection presents the inventory model at the in-plane spares level. It includes the demand model for in-plane spares, their resupply lead time, their backorders, and finally their mean stock level over a replenishment cycle.

\subsubsection{Demand model for in-plane spares: the satellite failures}
Satellite reliability is the factor at stake when designing a constellation maintenance strategy, as it is responsible for failures. In our approach, satellite failures are modeled by a Poisson distribution with the satellite failure rate as its parameter, meaning the number of failures per unit time \cite{collopy2003assigning}. The failure rate per constellation orbital plane is deduced from the satellite failure rate:

\begin{equation} \label{lplane}
\lambda_{\text{plane}} = \frac{N_{\text{sats}} \ \lambda_{\text{sat}} }{N_{\text{days}}}
\end{equation}

\noindent Note that an underlying assumption here is that the failed satellites are replaced by new spares immediately, which would be reasonable to assume for our optimal spare strategy.

\subsubsection{Resupply lead time from the parking orbits to the constellation orbits}\label{ltorbits}
As explained previously, the constellation's in-plane spare stocks get supplies from the closest available parking orbit at the time of the order, as the parking orbits drift relative to the constellation orbits. The lead time from the order processing by a constellation orbit to the actual delivery is therefore stochastic and its probabilistic distribution has to be derived. First, we need to determine the probability of a parking orbit to be available, and then derive the probability of a constellation orbit to get a supply from a specific parking orbit. Furthermore, the lead time distribution is derived from the geometry of the problem and orbital mechanics considerations.\\

\textit{(a) Probability of parking orbit availability}

The probability of parking orbit availability can be derived using a binomial-like distribution. The constellation orbits need to get a supply from the closest (i.e., minimum wait time) available parking orbit, while each parking orbit can either have available spare batches or be out-of-stock. Thus, given the probability of each parking orbit being available, $P_{\text{av}}$, we can derive the probability that a constellation orbit gets a supply from the $i^{th}$ closest parking orbit. 
Note that in our application, the probability that all parking orbits are out-of-stock at the time of delivery is very small, and thus can be neglected; therefore we assume that there is always one parking orbit available, which is not necessarily the closest one, to supply the in-plane stocks. 

The probability that a parking orbit has available spare batches, $P_{\text{av}}$, corresponds to the probability of visiting a parking orbit and not observing a stock-out. This probability is equal to the fraction of the demand not being backordered because the demand arrives at a constant rate. Therefore, $P_{\text{av}}$ can be expressed as follows: 

\begin{equation}
P_{\text{av}} = 1-\frac{ES_{\text{parking}}}{k_{\text{Q,parking}}}
\end{equation}

\noindent where $ES_{\text{parking}}$ is the expected number of backorders over a replenishment cycle at a parking orbit, which is derived in Section \ref{ESParking}. Note that since all of the parking orbits are analogous and evenly distributed, they are supposed to have the same $P_{\text{av}}$. 

Using this $P_{\text{av}}$, the probability of getting supply from the $i^{th}$ closest parking orbit is then obtained by summing all the possible cases:

\begin{equation}
\label{p_ith}
P(i^{th}) = \sum_{k = 1}^{N_{\text{parking}}-i+1} \ \binom{N_{\text{parking}} - i}{k-1} \ P_{\text{av}}^k \ (1-P_{\text{av}})^{N_{\text{parking}}-k}
\end{equation}

In order to demonstrate this expression, we consider a simple example. Assume that the chosen configuration is $N_{\text{parking}}=3$ and we want to determine the probability of getting supply from each parking orbit. 

\begin{itemize}
\item $1^{st}\  closest \ orbit$:  The possible cases and their respective probabilities are:
\begin{enumerate}
\item All orbits are available:
$ P = P_{\text{av}}^3 $
\item The $1^{st}$ closest orbit is available, the  $2^{nd}$ is available and the $3^{rd}$ is not available:
$P = P_{\text{av}}^2 \ (1-P_{\text{av}})$
\item The $1^{st}$ closest orbit is available, the  $2^{nd}$ is not available and the $3^{rd}$ is available:
$ P = P_{\text{av}}^2 \ (1-P_{\text{av}})$
\item The $1^{st}$ closest orbit is available, the  $2^{nd}$ and $3^{rd}$ orbits are not available:
$P = P_{\text{av}} \ (1-P_{\text{av}})^{2}$
\end{enumerate}

So $$P(1^{st}) = P_{\text{av}}^3  + 2  (P_{\text{av}}^2 \ (1-P_{\text{av}})) + P_{\text{av}} \ (1-P_{\text{av}})^{2} = \sum_{k = 1}^{3} \binom {3 - 1}{k-1} \ P_{\text{av}}^k \ (1-P_{\text{av}})^{3-k}$$

\item $2^{nd}\  closest \ orbit$:  The possible cases and their respective probabilities are:
\begin{enumerate}
\item The $1^{st}$ closest orbit is not available, the  $2^{nd}$ is available and the $3^{rd}$ is available:
$ P = P_{\text{av}}^2 \ (1-P_{\text{av}})$
\item The $1^{st}$ closest orbit is not available, the  $2^{nd}$ is available and the $3^{rd}$ is not available:
$ P = P_{\text{av}} \ (1-P_{\text{av}})^{2} $
\end{enumerate}

So $$P(2^{nd}) = P_{\text{av}}^2 \ (1-P_{\text{av}}) + P_{\text{av}} \ (1-P_{\text{av}})^{2} = \sum_{k = 1}^{2} \binom {3 - 2}{k-1} \ P_{\text{av}}^k \ (1-P_{\text{av}})^{3-k}$$

\item $3^{rd}\  closest \ orbit$:  The only possible case and its probability are: 
\begin{enumerate}
\item The $1^{st}$ and $2^{nd}$ closest orbits are not available and the $3^{rd}$ is available:
$ P = P_{\text{av}} \ (1-P_{\text{av}})^{2}$
\end{enumerate}

So $$P(3^{rd}) = P_{\text{av}} \ (1-P_{\text{av}})^{2} = \sum_{k = 1}^{1} \binom{3 - 3}{k-1} \ P_{\text{av}}^k \ (1-P_{\text{av}})^{3-k} $$
\end{itemize}

\textit{(b) Lead time distribution}

The spare model presented in this paper takes advantage of the RAAN drift caused by Earth's gravitational field (see Section \ref{Orbitpert}). Over time, a parking orbit will visit all the constellation orbits and hence is able to service failures in all of them. When a parking orbit and the constellation orbit of interest are aligned, the orbital maneuver becomes feasible and a transfer is performed (see Section \ref{Orbittransfer} for details about the transfer).
The lead time to transfer batches of satellites from the parking orbits to the constellation orbits is the result of the drift time to align the orbital planes and the actual time of flight \cite{palmade1998skybridge}.

\indent A probability distribution now has to be defined to describe the transfer time, meaning the lead time from the parking orbits to the constellation orbits. Spares are transferred from the closest parking orbit with available supply at the time of the order. As the parking orbits are angularly equally distributed, it divides the possible RAAN differences for drift into $N_{\text{parking}}$ intervals: $[0, \frac{2\pi}{N_{\text{parking}}}], [\frac{2\pi}{N_{\text{parking}}},\frac{4\pi}{N_{\text{parking}}}], ..., [\frac{2\pi(N_{\text{parking}}-1)}{N_{\text{parking}}},2\pi]$. Indeed, if spares are transferred from the closest parking orbit to the constellation orbit of interest, the possible RAAN differences ($\Delta \Omega$) belong to $[0, \frac{2\pi}{N_{\text{parking}}}]$, while if the spares are transferred from the $i^{th}$ closest parking, $\Delta \Omega \in [(i-1)\frac{2\pi}{N_{\text{parking}}},i\frac{2\pi}{N_{\text{parking}}}]$. Given that the drift rates are fixed by the semi-major axis and the inclination and that the parking orbits are equally distributed, we can consider that transfer times are uniformly distributed in each possible interval (see Eq. \ref{ltuniform}).

\begin{equation}
\label{ltuniform}
T_{\text{plane}}(i^{th}) \sim \mathcal{U}\left\{t_\text{transfer}(\Delta \Omega = (i-1)\frac{2\pi}{N_{\text{parking}}}), t_\text{transfer}(\Delta \Omega = i\frac{2\pi}{N_{\text{parking}}})\right\}
\end{equation}

\noindent where $t_\text{transfer}(\Delta \Omega)$  is the summation of the drift waiting time  for an angular difference of $\Delta \Omega$ and the time of flight, each of which can be calculated using Eqs. \ref{drift} and \ref{tof}, respectively. With $P(i^{th})$ found in Eq. \ref{p_ith} and $T_{\text{plane}}(i^{th}) $ found in Eq. \ref{ltuniform}, we can find the lead time distribution from the parking orbits to the constellation orbits.

\subsubsection{Backorders at the constellation orbits} \label{ES orbit}
For the in-plane spare stocks, the expected number of backorders over a cycle, $ES_{\text{plane}}$, can be calculated from the distribution of lead time and the expected demand during this lead time \cite{ganeshan1999managing}.
\begin{equation} \label{ESplane}
ES_{\text{plane}} = \int_{T_{\text{plane}}} \ ES_{ T_{\text{plane}}}(s_{\text{plane}}) \ f_{\text{plane}}(T_{\text{plane}}) \ dT_{\text{plane}}
\end{equation}

\noindent where $\ ES_{\tau}(s_{\text{plane}})$ is the expected backorders for the lead time being $\tau$ and the threshold stock level being $s_{\text{plane}}$, and $f_{\text{plane}}$ is the probability density function of the lead time to the constellation orbits found in Section \ref{ltorbits}. Since an $(s,Q)$ policy is considered, the expected backorders can be found using the approach in Section \ref{EStau}. With this $ES_{\text{plane}}$, we can find the order fill rate using Eq. \ref{fillrateplane}.

\begin{equation} \label{fillrateplane}
\rho_{\text{plane}} =  1-\frac{ES_{\text{plane}}}{Q_{\text{plane}}}
\end{equation}

\subsubsection{Mean stock level of in-plane spares}
Finally, the mean stock level of spare parts should be evaluated to further calculate the holding cost of the spare strategy. The resulting mean stock level of in-plane spares over a cycle is calculated as the expected mean stock level over all possible lead times. According to the theory in Section \ref{subsecstock}, this mean stock level is given by Eq. \ref{Splane}.
\begin{equation} \label{Splane}
\overline{Stock_{\text{plane}}}= \int_{T_{\text{plane}}} \left\{\frac{Q_{\text{plane}}}{2} + s_{\text{plane}} - N_{\text{fail,plane}}(T_{\text{plane}}) +\frac{1}{2}\right\} f_{\text{plane}}(T_{\text{plane}})  \ dT_{\text{plane}}
\end{equation}

\noindent Note that even though the cycle length is stochastic, this mean stock level over a given cycle is equal to the mean stock level over the entire time horizon because the cycle length distribution (governed by the demand rate) is independent of the lead time distribution. Also, an underlying assumption is that the backorders are negligible, which is a reasonable assumption for our optimal spare strategy. 

\subsection{Parking Spares Inventory Model}\label{parkingmodel}
The inventory model at the parking orbits also follows an $(s,Q)$ policy. This subsection presents the inventory model at the parking spares level. It includes the demand model for parking spares, their resupply lead time, their backorders, and finally their mean stock level.

\subsubsection{Demand model for parking spares}\label{parkingdemand}
The demand process at the spare parking orbits is derived from the failure process and policy model at the constellation orbits.
Looking at the ordering process from one constellation orbital plane, an order is placed every $Q_{\text{plane}}$ failures on average and those failures are Poisson distributed. Therefore, the times between consecutive orders from this plane are Erlang-$Q_{\text{plane}}$ distributed according to the relationship between the two stochastic distributions \cite{ganeshan1999managing}. The orders placed at all spare parking orbits combined is the superposition of the orders from all constellation orbits. When $N_{\text{plane}}$ is sufficiently large (meaning $N_{\text{plane}} \geq 20$), the superposition of those $N_{\text{plane}}$ Poisson processes can also be considered as a Poisson process, with rate $N_{\text{plane}}\frac{\lambda_{\text{plane}}}{Q_{\text{plane}}}$ \cite{ccinlar1975exceptional}\cite{schwarz1985fill}. Considering the symmetry of the problem where all spare parking orbits are equally distributed, each parking orbit is thus subject to a Poisson demand with rate $\lambda_{\text{parking}}$, derived in Eq. \ref{lpark}.

\begin{equation}
\label{lpark}
\lambda_{\text{parking}} = N_{\text{plane}} \ \frac{\lambda_{\text{plane}}}{Q_{\text{plane}}} \ \frac{1}{N_{\text{parking}}}
\end{equation}

\subsubsection{Resupply lead time from the ground to the parking orbits}\label{ltparking}
The spare parking orbits are replenished from the ground using rocket launches, with a certain lead time denoted as $T_{\text{parking}}$. This lead time takes into account the launch order processing time and the waiting time for the next launch window. The model proposed in this paper does not include any manufacturing delay, assuming the spare stock on the ground to be always available. The order processing time is considered to be constant, while the waiting time for the next launch window is assumed to be exponentially distributed in accordance with launch schedules databases (see Appendix A).

\begin{equation}
T_{\text{parking}} \sim \mathcal{E}(\mu_{\text{launch}}) + pt_{\text{launch}}
\end{equation}
\noindent where $ \mathcal{E}(\mu_{\text{launch}})$ is the exponential distribution with mean $\mu_{\text{launch}}$.

\subsubsection{Backorders at the parking orbits}\label{ESParking}
The inventory policy for the parking spare stocks is similar to the one used for the in-plane spare stocks. Therefore, the expected number of backorders over a replenishment cycle at a parking orbit in units of batches $Q_{\text{plane}}$ can be derived using the same technique as used in Section \ref{ES orbit}, and is given by Eq. \ref{ESparking}:
\begin{equation} \label{ESparking}
ES_{\text{parking}} = \int_{T_{\text{parking}}} \ ES_{T_{\text{parking}}}(k_{\text{s,parking}}) \ f_{\text{parking}}(T_{\text{parking}}) \ dT_{\text{parking}}
\end{equation}
With this $ES_{\text{parking}}$, we can find the order fill rate using Eq. \ref{fillrateparking}.

\begin{equation} \label{fillrateparking}
\rho_{\text{parking}}=1-\frac{ES_{\text{parking}}}{k_{\text{Q,parking}}}
\end{equation}

\subsubsection{Mean stock level of parking spares}
The replenishment cycles at the parking orbits follows the same characteristics as the in-plane spares cycles. Therefore, the mean stock level at the parking orbits is, in units of batches $Q_{\text{plane}}$:
\begin{equation} \label{Sparking}
\overline{Stock_{\text{parking}}} = \int_{T_{\text{parking}}} \left\{\frac{k_{\text{Q,parking}}}{2} + k_{\text{s,parking}} - N_{\text{fail,parking}}(T_{\text{parking}})+\frac{1}{2}\right\} \\ f_{\text{parking}}(T_{\text{parking}})  \ dT_{\text{parking}}
\end{equation}

\noindent where $N_{\text{fail,parking}}(T_{\text{parking}})$ is the failure demand at the parking orbits over the lead time $T_{\text{parking}}$, in units of batches $Q_{\text{plane}}$.

\subsection{Total Cost Model} \label{totcostmodel}
The goal of the model is to estimate the cost of the spare strategy to maintain the system. For this purpose, four types of costs are considered: the manufacturing ($c_{\text{manufacturing}}$), holding($c_{\text{holding}}$), launching ($c_{\text{lauch}}$), and orbital maneuvering ($c_{\text{maneuvering}}$) costs. The sum of the aforementioned cost items gives us the total expected spare strategy annual cost (TESSAC):

\begin{equation} \label{TESSAC}
TESSAC = c_{\text{manufacturing}}+c_{\text{holding}}+c_{\text{launch}}+c_{\text{maneuvering}}
\end{equation}

\subsubsection{Manufacturing cost}
The annual manufacturing cost is derived from the total number of failures observed over a year. As the failures are Poisson distributed with a rate $\lambda_{\text{plane}}$ for each of the $N_{\text{plane}}$ planes, $c_{\text{manufacturing}}$ is given by:

\begin{equation}
c_{\text{manufacturing}} = p_{\text{sat}} \ \lambda_{\text{plane}} \ N_{\text{plane}} \ N_{\text{days}}
\end{equation}
\noindent where  $\lambda_{\text{plane}}$  is derived in Eq. \ref{lplane}.

\subsubsection{Holding cost}
The annual holding cost is associated with the spare strategy. Having spare satellites in orbits represents a substantial cost due to their operations and station keeping. The annual holding cost of in-space and parking spare satellites is defined using the mean spare stock level at each orbit. 

\begin{equation}
c_{\text{holding}} = p_{\text{holding}} \ \{\overline{Stock_{\text{plane}}} \ N_{\text{plane}} + \overline{Stock_{\text{parking}}} \ Q_{\text{plane}} \ N_{\text{parking}} \} 
\end{equation}
\noindent where $ \overline{Stock_{\text{plane}}}$ and $ \overline{Stock_{\text{parking}}}$ are given by Eq. \ref{Splane} and Eq. \ref{Sparking}, respectively.

\subsubsection{Launch cost}
The annual launch cost is derived from the demand generated at the parking orbits:

\begin{equation}
\label{launcheq}
c_{\text{launch}} = p_{\text{launch}} \frac{\lambda_{\text{parking}}\ Q_{\text{plane}} }{Q_{\text{parking}}} \ N_{\text{parking}} \ N_{\text{days}}
\end{equation}

\noindent where $Q_{\text{parking}}$ is given by Eq. \ref{batcheq}, $\lambda_{\text{parking}}$ is given by Eq. \ref{lpark}, and $p_{\text{launch}}$ is the launch cost given by Eq. \ref{minlaunch}. Two possibilities are offered regarding the launch of spare satellites, mimicking the launch options considered by OneWeb \cite{onewebLaunch} : 
\begin{enumerate}
\item Using a full capacity rocket, allowing launches up to the rocket capacity, $cap_{\text{launch}}$ satellites, at once for a fixed cost $p_{\text{launch,full}}$, which does not depend on the actual batch number of satellites effectively launched from this rocket.
\item Using a unit-satellite launcher at a cost of $p_{\text{launch,unit}}$ per launch, i.e., per spare satellite launched. Given the specificity of this type of launcher, which is not dependent on government maintained launch ranges to launch, it is considered possible to launch several rockets at the same time \cite{virgin}. This option requires as many launchers as the number of satellites that need to be launched.
\end{enumerate}
\begin{equation}
\label{minlaunch}
p_{\text{launch}} = \min \left\{p_{\text{launch,full}}, \ Q_{\text{parking}} \ p_{\text{launch,unit}}\right\}
\end{equation}

\subsubsection{Maneuvering cost}
The annual maneuvering cost corresponds to the fuel mass required to perform maneuvers for all orbital transfers required over a year, affected by a conversion coefficient $\epsilon_{\text{maneuvering}} \ \mathrm{ [million \ US\$/kg]}$.

\begin{equation}
c_ {\text{maneuvering}} =  m_{\text{fuel}}\ \lambda_{\text{plane}} \ N_{\text{plane}} \ N_{\text{days}}  \ \epsilon_{\text{maneuvering}}
\end{equation}
\noindent where $m_{\text{fuel}}$ is calculated in  Eq. \ref{mfuel}  and $\lambda_{\text{plane}}$ is given by Eq. \ref{lplane}.

\section{Assessment of Model Accuracy}\label{validationsection}
The model presented in the previous section is an analytical model which allows computationally efficient evaluation of a spare policy, even for mega-scale constellations. Nevertheless, it relies on a number of simplifying assumptions (e.g., demand distribution at the parking orbits, low probability of backorders), and its accuracy needs to be assessed using simulations. Those simulations are performed with a variety of values for parameters and variables. Once the accuracy of the model is shown to be acceptable, it can be used for optimization of the spare policy without relying any more on costly simulations.

A set of 25 unique test problems are constructed using Latin Hypercube Sampling (LHS). This method generates near-random sets of parameters from a multidimensional trade space, forcing the samples to represent the real variability of the parameters \cite{stein1987large}. The parameters used in all the simulation experiments are given in Table \ref{simparams}. They are representative of mega-constellation figures such as OneWeb \cite{onewebBuilt}. The sampled trade space can be found in Table \ref{lhs}. The reorder points $s_{\text{plane}}$ and $k_{\text{s,parking}}$ for simulations are determined through the analytical model for a set of requirements on the order fill rates:  
\begin{equation}
	\rho_{\text{plane}}^{N_{\text{plane}}}\geq 0.95
	\end{equation}
	\begin{equation} \rho_{\text{parking}}^{N_{\text{parking}}} \geq 0.95
\end{equation}
\noindent where $\rho_{\text{plane}}$ is calculated in Eq. \ref{fillrateplane} and $\rho_{\text{parking}}$ is calculated in Eq. \ref{fillrateparking}.
These requirements are set because we are only interested in the highly efficient policies with few backorders, and that is also an underlying assumption for the analytical model. The results from the simulations using these $(s,Q)$ policies are used to assess the accuracy of outputs from the analytical model developed in Section \ref{modelsection}: the mean stock level of in-plane spares, the mean stock level of parking spares, the order fill rate of the in-plane spare stocks, the order fill rate of the parking spare stocks, and the TESSAC. Each simulation is run for 15 years and each case includes 100 simulations. Given the simulation and modeling results, relative errors percentages are calculated according to:
  
  \begin{equation}
  \label{error}
  \frac{|value_{\text{sim}} - value_\text{model}|}{value_\text{sim}} *100
  \end{equation}

  \begin{table}[ht!]
  	
  	\centering
  	\renewcommand{\arraystretch}{1.3}
  	\caption{Fixed simulation parameters}
  	\begin{tabular}{cccc}
  		\hline
  		\hline
  		Parameter & Notation & Value & Unit  \\
  		\hline
  		Fuel mass conversion coefficient & $\epsilon_{\text{maneuvering}}$ & 0.001 & million US\$/kg\\
  		Annual satellite holding cost & $p_{\text{holding}}$ & 0.5& million US\$/satellite/year \\
  		Launch capacity& $cap_{\text{launch}}$ & 34 & satellites \\
  		Satellite dry mass &  $m_{\text{dry}}$& 150 & kg \\
  		Satellite manufacturing cost & $p_{\text{sat}}$ & 0.5 & million US\$/unit  \\
  		Full rocket launch price & $p_{\text{launch,full}}$ & 47.6 & million US\$/launch  \\
  		Unique satellite rocket launch cost & $p_{\text{launch,unit}}$ & 10 & million US\$/launch  \\
  		Effective exhaust velocity & $v_\text{exhaust}$ & 2.16 & km/s  \\
  		\hline
  		\hline
  	\end{tabular}
  	\label{simparams}
  \end{table}

  \begin{table}[ht!]
  	
  	\centering
  	\renewcommand{\arraystretch}{1.3}
  	\caption{Sampled trade space for LHS}
  	\begin{tabular}{cccc}
  		\hline
  		\hline
  		Parameter& Notation & Bounds & Unit \\
  		\hline
  		Launch order processing time & $pt_{\text{launch}}$ & $30 \le pt_{\text{launch}} \le 120 $ & days  \\
  		Constellation orbit altitude & $h_{\text{plane}}$ & $1000 \le h_{\text{plane}} \le 2000$  &  km \\
  		Parking orbit altitude & $h_{\text{parking}}$ & $700 \le h_{\text{parking}} \leq 1000$  &  km \\
  		Inclination &  $i$ & $30 \le i \le 70$ & deg  \\
  		Satellite failure rate &  $\lambda_{\text{sat}}$ & $0.001 \le \lambda_{\text{sat}} \le 0.1$ & failures/year  \\
  		Mean time between launches & $\mu_{\text{launch}}$ & $30 \le \mu_{\text{launch}} \le 90$ & days  \\
  		Number of planes in the constellation & $N_{\text{plane}}$ & $20 \le N_{\text{plane}} \le 40$ & planes  \\
  		Number of parking orbits & $N_{\text{parking}}$ & $1 \le N_{\text{parking}} \le 20$ & planes  \\
  		Number of operational satellites per orbital plane &     $N_{\text{sats}}$ & $20 \le N_{\text{sats}} \le 60$ & satellites/plane  \\             
  		Order batch size for in-plane spares & $Q_{\text{plane}}$ & $1 \le Q_{\text{plane}} \le 10$ & satellites  \\
  		Order batch size for parking spares &   $k_{\text{Q,parking}}$ & $1 \le k_{\text{Q,parking}} \le 10$ & $Q_{\text{plane}}$  \\
  		\hline
  		\hline
  	\end{tabular}
  	\label{lhs}
  \end{table}

The evaluation of the model through the relative percentage errors with simulations can be found in Table \ref{simerrors}. The results of the simulations indicate that the analytical model performs well, with relative error percentages ranging from 0.4\% to 4.1\% on average. 
The mean stock levels of in-plane spares and parking spares reveal relative errors of 1.7\% and 4.1\%, respectively. Those low error percentages indicate that the model accurately estimates the stocks given the lead time distributions. 
The order fill rates of the in-plane spare stocks and the parking spare stocks are very well estimated by the model with relative errors of 0.8\% and 0.4\%, respectively. The calculation of the expected backorders of replenishment cycles through demand and lead time distributions is therefore accurately performed by the analytical model. 
Finally, the TESSAC error is quantified, leading to a relative error of 1.6\% on average.  
These results indicate that the developed analytical model can approximate the simulated values well without running computationally costly simulations. \footnote{A typical set of 100 simulations over 15 years takes more than one hour, whereas the analytical model takes less than ten seconds, both with MATLAB R2016a on an Intel Core i5-6300U 2.4 GHz platform.}

  \begin{table}[ht!]
  	\centering
  	
  	\renewcommand{\arraystretch}{1.3}
  	\caption{Averaged relative errors percentages of the analytical model vs. simulations}
  	\begin{tabular}{cc}
  		\hline
  		\hline
  		Output & Relative error percentage\\
  		\hline
  		Mean stock level of in-plane spares & 1.7\% \\
  		Mean stock level of parking spares & 4.1\%\\
  		Order fill rate of the in-plane spare stocks  & 0.8\% \\
  		Order fill rate of the parking spare stocks  & 0.4\% \\
  		TESSAC  & 1.6\% \\
  		\hline
  		\hline
  	\end{tabular}
  	\label{simerrors}
  \end{table}
  
\section{Optimization Problem Formulation}\label{optimizationsection}
With the developed model, we can develop an optimization problem formulation to find the optimal spare strategy. The spare strategy design problem can be formulated as a mixed-integer nonlinear problem comprising six variables. The objective of the optimization problem is to design a spare strategy which minimizes the TESSAC for a given operational constellation.

\subsection{Variables}

Table \ref{variables} presents the spare strategy decision variables along with their possible ranges of values and integer constraints. 

\begin{table}[ht!]
\centering

\renewcommand{\arraystretch}{1.3}
  \caption{Optimization design variables}
\begin{tabular}{cccc}
  \hline
  \hline
  Variable & Unit & Bounds & Constraint \\
  \hline
  $N_{\text{parking}}$  & - & $1 \leq N_{\text{parking}} \leq 20$ & integer \\
  $h_{\text{parking}}$ &  km &  $700 \leq h_{\text{parking}} \leq 1000$ & -\\
  $Q_{\text{plane}}$ &  satellites& $1 \leq Q_{\text{plane}} \leq 10$ & integer \\
  $s_{\text{plane}}$ &  satellites & $1 \leq s_{\text{plane}} \leq 10$ & integer \\
  $k_{\text{Q,parking}}$ &  $Q_{\text{plane}}$ & $1 \leq k_{\text{Q,parking}} \leq 10$ & integer\\
  $k_{\text{s,parking}}$ &  $Q_{\text{plane}}$ & $1 \leq k_{\text{s,parking}} \leq 10$ &  integer \\
    \hline
    \hline
\end{tabular}
    \label{variables}
\end{table}

From the specific formulation of our problem, it is important to note two major implications of the parking orbit design choice:

 \begin{enumerate}
 \item The number of spare parking orbits $N_{\text{parking}}$ determines the maximum angular difference observed between the parking orbits and the constellation orbits. While a large number of parking orbits results in shorter transfer times, it can also lead to higher costs. 
 \item The altitude of the spare parking orbits $h_{\text{parking}}$ determines the relative rotation of the two orbits and, consequently, the drift time required to carry out the actual transfer of spares from the parking orbits to the constellation orbits. It also, to a smaller extent, influences the time of flight of the maneuver. 
\end{enumerate}

\subsection{Objective Function}
The optimization of the spare strategy is made to minimize the TESSAC, comprising the costs of manufacturing, holding, launching, and orbital maneuvering of the spare satellites over a year of maintenance:
\begin{equation}
min_{x = [N_{\text{parking}},h_{\text{parking}},Q_{\text{plane}},s_{\text{plane}},k_{\text{Q,parking}},k_{\text{s,parking}}]} \ J(x) \ = \ TESSAC(x)
\end{equation}

\noindent where $TESSAC$ is given by Eq. \ref{TESSAC} according to the analytical model detailed in Section \ref{modelsection}.

\subsection{Constraints} 
The constraints for the optimization problem have two components. 

The first component is to enforce the launch capacity constraint. Since every launch vehicle can only deliver up to $cap_{\text{launch}}$ satellites, we have the following constraint:
\begin{equation}
Q_{\text{parking}} \le cap_{\text{launch}}
\label{capacity}
\end{equation}
\noindent where $Q_{\text{parking}}$ is a function of $k_{\text{Q,parking}}$ and $Q_{\text{plane}}$ according to Eq. \ref{batcheq}.

The second component is to ensure the multi-echelon spare policy to meet a global requirement for efficiency $\rho_T$. This global efficiency can be achieved using different relative configurations between in-plane spares and parking orbit spares, thus allowing more flexibility in the design of the inventory model at different echelons. The constraints can be written as follows:  
 \begin{equation}
\rho_T \le  \rho_{\text{plane}}^{N_{\text{plane}}} \ \rho_{\text{parking}}^{N_{\text{parking}}} 
 \end{equation}
\noindent where $\rho_{\text{plane}}$ is calculated in Eq. \ref{fillrateplane} and $\rho_{\text{parking}}$ is calculated in Eq. \ref{fillrateparking}. This constraint limits the backorders, making them negligible for the mean stock level calculation as described in Section \ref{secinventory}.

\subsection{Optimizer}
The optimization has to be performed using a mixed-integer nonlinear solver to meet the formulation of the problem. For the purpose of this paper, the single objective genetic algorithm (GA) embedded in Matlab is used to complete the optimization.

\section{Numerical Example}\label{examplessection}
This section shows a numerical example of satellite mega-constellation spare strategy optimization. Specifically, we focus on evaluating the value of parking orbits utilizing our proposed inventory model. The specific strategy of using parking orbits drifting and supplying the constellation orbits have been proposed in the existing literature; however, no study has been able to optimize the operational strategy of these parking orbits in a scalable and rigorous way. Thus, it is of interest to evaluate the benefits of having parking orbits in our spare strategy design. A competitive design comprising only in-plane spares replenished directly from ground rocket launches with no parking orbits is also optimized for an $(s,Q)$ policy, given the same parameters and satellite configuration. Note that, for the in-plane-only strategy, the upper bound for $Q_{\text{plane}}$ is specified by the launch capacity constraint because each rocket only delivers one batch to a constellation orbit. A cost comparison is established between the in-plane-only and  multi-echelon strategies.

\subsection{Mega-Constellation Configuration and Requirements}
The implementation of a study case for a LEO satellite mega-constellation is described, for which an optimization using the model previously exposed is performed. Given the nominal constellation configuration and performance requirements, the optimizer derives the best set of variables $[N_{\text{parking}}, h_{\text{parking}},Q_{\text{plane}},s_{\text{plane}},\allowbreak k_{\text{Q,parking}}, k_{\text{s,parking}}]$ with respect to the objective fitness function $J$. The used parameters remain the same as in Table \ref{simparams} and the chosen LEO configuration and performance requirements are: 
\begin{equation*}
\left\{
\begin{array}{cc}
  h_{\text{plane}} = 1200 \ \mathrm{km} \\
	i = 50 ^{\mathrm{o}}   \\
    N_{\text{plane}} = 40  \ \mathrm{planes} \\
    N_{\text{sats}} = 40 \ \mathrm{satellites/plane} \\
    \lambda_{\text{sat}} = 0.05 \ \mathrm{failures/year}\\
    \rho_{T} = 0.95
\end{array}
\right.
\end{equation*}
Specific parameters related to launch are:
\begin{equation*}
\left\{
\begin{array}{cc}
  \mu_{\text{launch}} = 66.7 \ \mathrm{days} \\
	pt_{\text{launch}} = 90 \ \mathrm{days}   \\
\end{array}
\right.
\end{equation*}

\subsection{Results and Analysis}

The results of the optimization for both in-plane-only and  multi-echelon strategies are summarized in Table \ref{OptimSummary}, along with a comparison of their respective TESSAC.

\begin{table}[ht!]
\centering
\renewcommand{\arraystretch}{1.3}

\caption{Optimization results and comparison between the  in-plane-only and multi-echelon spare strategies}
\begin{tabular}{ccccc}
\hline
\hline
Strategy & Chromosome & Optimal solution &  \begin{tabular}[x]{@{}c@{}}TESSAC \\ (million US\$/year)\end{tabular}& Comparison \\
\hline
\begin{tabular}[x]{@{}c@{}}In-plane-only \\ (traditional)\end{tabular} &  $[Q_{\text{plane}},s_{\text{plane}}]$ & $[20,\ 4]$ & $503.2$ & - \\
Multi-echelon &  \begin{tabular}[x]{@{}c@{}}$[N_{\text{parking}},h_{\text{parking}},$ \\ $Q_{\text{plane}},s_{\text{plane}},k_{\text{Q,parking}},k_{\text{s,parking}}]$\end{tabular}& \begin{tabular}[x]{@{}c@{}}$[3,\ 792.3,\ 4,\ 3,\ 8,\ 8]$\end{tabular} & $319.1$ & $-36.6\%$ \\
\hline
\hline
\end{tabular}
\label{OptimSummary}
\end{table}

\indent The chosen design for the multi-echelon spare strategy has three parking orbits at an altitude of 792.3 km with the $(s_{\text{parking}},Q_{\text{parking}}) = (32,32)$ inventory policy, along with the $(s_{\text{plane}},Q_{\text{plane}}) = (3,4)$ policy at each orbit's spare stock. The associated TESSAC is $ J^*(x^*)= \mathrm{US\$} \  319.1 \ \mathrm{million/year}$. 
In comparison to this chosen design, the in-plane-only strategy optimization leads to an inventory policy of $(s_{\text{plane}},Q_{\text{plane}})= (4,20)$ and the associated annual maintenance cost is $J^*(x^*_\text{in-plane-only})= \mathrm{US\$} \ 503.2 \ \mathrm{ million/year}$.

The results of the performed optimization show interesting features.

\indent First, the comparison of the multi-level mixed strategy with a single in-plane-only strategy shows the value of introducing another level of constellation spares as parking orbits and optimally designing its inventory management, reducing the TESSAC by 36.6\%. 

Furthermore, parking orbits allow us to take full advantage of the batch launch discount effectively, which is captured thanks to our unique optimization framework. Indeed, spare satellites can be launched in large quantities to the parking orbits as they will supply all the constellation orbits, whose demand rate is high. On the contrary, if large batches of spare satellites are launched directly to a constellation orbital plane, they will service only in-plane failures for that specific orbital plane, whose demand rate is much lower than that of the parking orbits. As a result, given a similar batch launch quantity, launching directly to the in-plane stocks (i.e., not having parking orbits) can result in higher costs primarily due to the associated holding costs.

Related to the above point, it is also worth noting that the optimal solution prefers a parking order quantity $Q_{\text{parking}}$ as close as possible to the launch capacity $cap_{\text{launch}}$. In fact, this parameter is set to $cap_{\text{launch}}=34$ and the results return $Q_{\text{parking}}=32$. Therefore, this $cap_{\text{launch}}$ plays a very important role in the search for the lowest possible maintenance policy and verifies the need to use satellite batch launches to reduce the cost of replenishment. Note that, although in this case $Q_{\text{parking}}$ almost matches $cap_{\text{launch}}$, this is a result of a tradeoff between the batch launch discount and the holding cost; it is expected that when the failure rate is very low, the optimizer would prefer to have less $Q_{\text{parking}}$ to save on holding cost.

Finally, the chosen parking orbit design $(N_{\text{parking}},h_{\text{parking}})$ also proves the value of having multiple parking orbits. The results show that the preferred design has three parking orbits. Indeed, even though having multiple parking orbits increases the costs of holding spare satellites, it also reduces the lead time to the constellation orbits; thus a sweet spot based on this tradeoff is found by the optimizer. Also, the altitude of the parking orbits (792.3 km) shows the compromise chosen by the optimizer between the duration of the lead time (especially the drift time to align the parking orbits and the constellation orbits) and the maneuver cost in terms of fuel mass required to perform the transfers. This demonstrates how our optimization can provide a direct impact on the design of a satellite constellation and its parking orbits.

\subsection{Sensitivity Analysis}

The key parameter for the analyzed constellation spare strategy optimization is the failure rate. In order to observe the effects of the failure rate on the optimized spare strategy solutions, a sensitivity analysis is performed for several values of failure rates.
As can be derived from \cite{dubos2010statistical} and \cite{erlank2015multicellular}, failure rates can range from $0.001$ to $0.9$ failures per year depending on the size of the spacecraft. Satellite constellations such as OneWeb and Starlink from SpaceX would fit in the "mini-satellite" category and thus displaying a failure rate of about $0.05$ failures per year after the first year.

The relative percentage of savings when using our unique multi-echelon approach using parking orbits compared to a single level of in-plane spares only is analyzed with respect to the TESSAC of each strategy. Fig. \ref{savings} shows the trend observed in savings with respect to different failure rates:

\begin{figure}[ht!]
\centering
\includegraphics[width = 13cm]{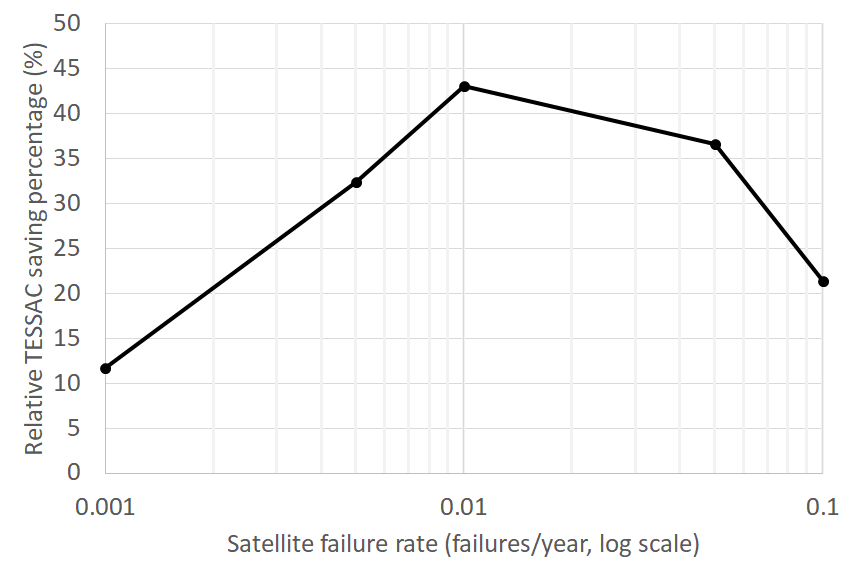}
\caption{Sensitivity analysis of the TESSAC savings using the multi-echelon strategy vs. the in-plane-spares-only strategy for different failure rates}
\label{savings}
\end{figure}

This sensitivity analysis shows that, for all cases, we observe cost savings when using the multi-echelon strategy as spares are better distributed and thus provide flexibility in the supply chain. In fact, spare satellites located in the multiple parking orbits are able to service all the constellation orbits and thus launched satellites are used more efficiently. The flow of spare satellites is more fluid as they do not get stuck in a particular plane, waiting for the next failure in this specific plane only. Even though the relative percentage of savings varies with the failure rate, the multi-echelon strategy is always preferred for the cases we tested.

Another interesting observation from Fig. \ref{savings} is that the largest cost saving is observed for the case with a medium failure rate ($\lambda_{\text{sat}}={0.01}\ \mathrm{failures/year}$), and the cases with higher or lower failure rates do not show as much cost savings by having parking orbits. This result can be interpreted as follows:
\begin{itemize}
	\item When the failure rate is low, the spare demands are also small, and so the optimizer does not take advantage of the batch launch discount very often. Thus, the optimized multi-echelon strategy value less the option of having parking orbits (only one parking orbit is chosen for cases $\lambda_{\text{sat}}={0.001}\ \mathrm{failures/year}$). 
	As a result, the relative saving using the multi-echelon strategy is also relatively small. 
	\item In the case with a high failure rate, both strategies take advantage of the batch launch discount because of the large spare demands. Even if there are no parking orbits, satellites could still be launched in batches directly to the constellation orbits to satisfy the demands. This configuration provides a relatively small saving using the multi-echelon strategy.
	\item The largest saving is observed for medium failure rates. The multi-echelon strategy takes full advantage of the batch launch discount, whereas in-plane-only strategy does not. The benefit of having parking orbits is the largest in this case, up to approximately 43\% of cost saving.
\end{itemize}

\section{Conclusion}\label{conclusionsection}
This paper presented a novel model for satellite constellation spare strategies using a multi-echelon inventory approach, and proposed an optimization formulation using this model to minimize the total annual cost of the spare strategy policy. The model views the satellite constellation spare strategy as a multi-level spare supply chain system, comprising the ground (supplier), multiple parking orbits (warehouses), and multiple orbital planes (retailers), all ruled by (s,Q) inventory policies and under the assumption of stochastic demand (failures) and lead times. 
Our inventory model is unique in that it has multiple drifting warehouses (parking orbit spare stocks), which are all capable of resupplying all the retailers (in-plane spare stocks), and the actual resupply pathway is chosen according to the availability and the lead time distribution. 
The measures of performance for a chosen spare strategy are derived from the analytical model, and a cost model of a strategy is developed, including manufacturing, holding, launch, and maneuvering costs. The accuracy of the proposed model is assessed using simulations. The paper additionally developed an optimization problem formulation to minimize the cost of maintenance under performance requirements, and the numerical case study demonstrated the value of having this multi-echelon mixed-strategy spare strategy for satellite mega-constellations. The importance of the batch launch discount is stressed in those results, along with the flexibility conveyed by the multiple parking orbits being able to deliver spares to all orbital planes. 

This research can be further extended in multiple directions. First, using non-identical parking orbits and non-identical orbital planes policies could allow more flexibility in the system to provide the same required efficiency. Also, the model presented in this paper assumes that ground spares are always available to launch with a given lead time, which is a reasonable assumption given the current satellite production rates; however, the possibility of ground spares to be out-of-stock could also be incorporated for a more accurate representation of reality. Furthermore, although this paper considered a hard efficiency constraint, an alternative approach would be to include in the objective function the cost of different efficiency outcomes so that the model itself can make the tradeoff on the optimal level of efficiency. Finally, this paper supposes a Poisson distribution of failures; however, existing satellite reliability analysis exhibited the problem of infant mortality \cite{castet2009satellite}
and introduced the use of "degraded states" \cite{castet2010beyond}. Therefore, more realistic consideration of satellite failure could be implemented using those observations.

\appendix

\section*{Appendix A: Launch time distribution}
\renewcommand{\thefigure}{A\arabic{figure}}
\setcounter{figure}{0}
Based on the launch data retrieved from \citep{soyouz} and \citep{falcon}, an exponential distribution fit is derived for the times between two consecutive successful launches. The example of the Soyuz rocket launches is given in Fig. \ref{launchdata}, where the obtained exponential parameter (i.e., the average time between two successful Soyuz launches) is 66.7 days.

\begin{figure}[ht!]
\centering
\includegraphics[width=18cm]{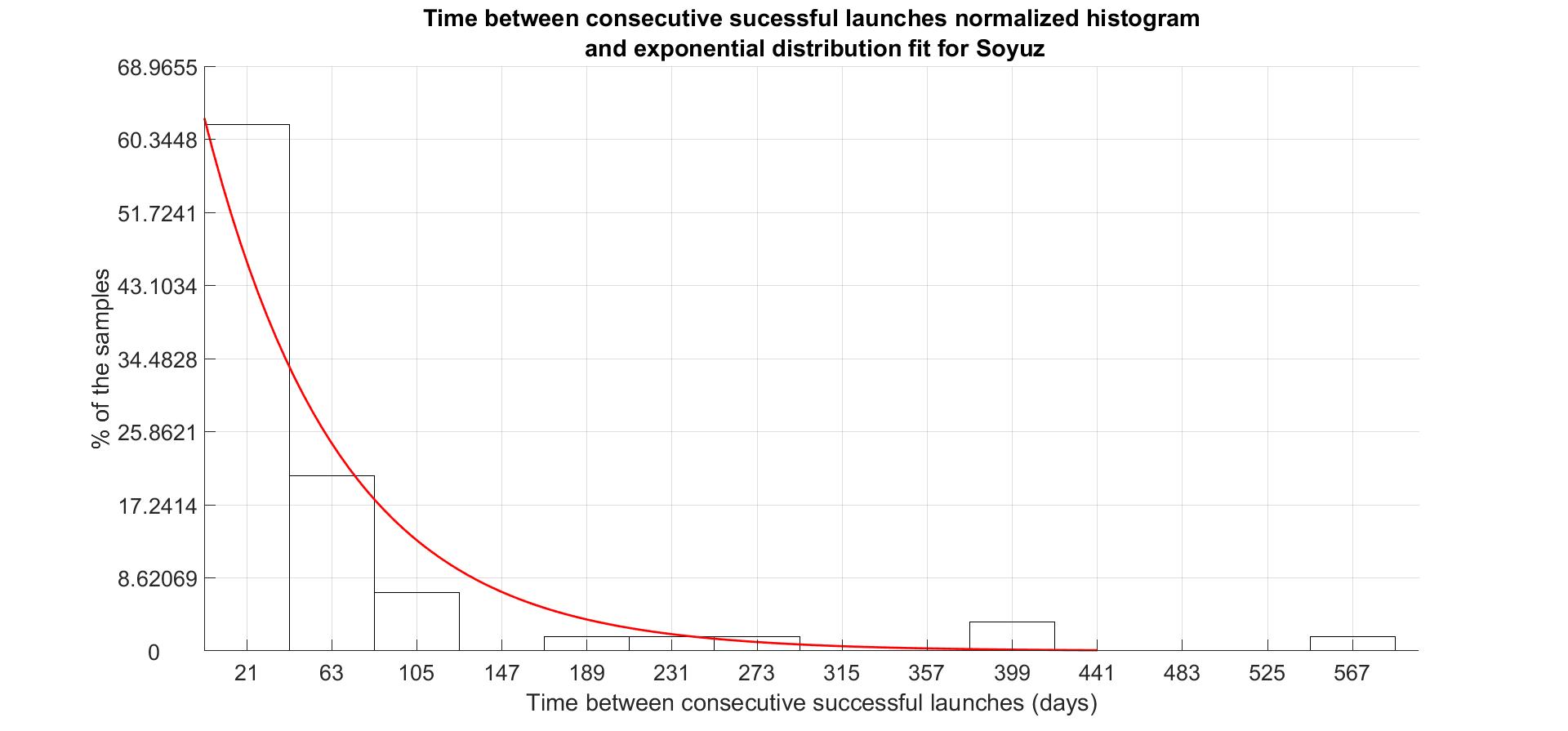}
\caption{Exponential distribution fit for Soyouz launches based on data from \cite{soyouz}}
\label{launchdata}
\end{figure}

\section*{Acknowledgment}
This research is supported by the Advanced Technology R\&D Center at Mitsubishi Electric Corporation. We appreciate Hang Woon Lee, Tiger Hou, Hao Chen, and Hai Wang for their reviews and thoughtful suggestions for improvement.

 \bibliography{AIAA_Spare}

\begin{thebibliography}{39}
\newcommand{\enquote}[1]{``#1''}
\providecommand{\natexlab}[1]{#1}
\providecommand{\url}[1]{\texttt{#1}}
\providecommand{\urlprefix}{URL }
\expandafter\ifx\csname urlstyle\endcsname\relax
  \providecommand{\doi}[1]{doi:\discretionary{}{}{}#1}\else
  \providecommand{\doi}{doi:\discretionary{}{}{}\begingroup
  \urlstyle{rm}\Url}\fi

\bibitem[{Ferringer and Spencer(2006)}]{ferringer2006satellite}
Ferringer, M.~P., and Spencer, D.~B., \enquote{Satellite constellation design
  tradeoffs using multiple-objective evolutionary computation,} \emph{Journal
  of Spacecraft and Rockets}, Vol.~43, No.~6, 2006, pp. 1404--1411.
\newblock Doi: \url{https://doi.org/10.2514/1.18788}.

\bibitem[{Ferringer et~al.(2007)Ferringer, Clifton, and
  Thompson}]{ferringer2007efficient}
Ferringer, M.~P., Clifton, R.~S., and Thompson, T.~G., \enquote{Efficient and
  accurate evolutionary multi-objective optimization paradigms for satellite
  constellation design,} \emph{Journal of Spacecraft and Rockets}, Vol.~44,
  No.~3, 2007, pp. 682--691.
\newblock Doi: \url{https://doi.org/10.2514/1.26747}.

\bibitem[{Budianto and Olds(2004)}]{budianto2004design}
Budianto, I.~A., and Olds, J.~R., \enquote{Design and deployment of a satellite
  constellation using collaborative optimization,} \emph{Journal of Spacecraft
  and Rockets}, Vol.~41, No.~6, 2004, pp. 956--963.
\newblock Doi: \url{https://doi.org/10.2514/1.14254}.

\bibitem[{Bandyopadhyay et~al.(2016)Bandyopadhyay, Foust, Subramanian, Chung,
  and Hadaegh}]{bandyopadhyay2016review}
Bandyopadhyay, S., Foust, R., Subramanian, G.~P., Chung, S.-J., and Hadaegh,
  F.~Y., \enquote{Review of formation flying and constellation missions using
  nanosatellites,} \emph{Journal of Spacecraft and Rockets}, Vol.~53, No.~3,
  2016, pp. 567--578.
\newblock Doi: \url{https://doi.org/10.2514/1.A33291}.

\bibitem[{Airbus({})}]{OneWebPic}
Airbus, \enquote{OneWeb satellites,}
  \url{http://www.airbus.com/space/telecommunications-satellites/oneweb-satellites-connection-for-people-all-over-the-globe.html},
  {}.
\newblock Accessed 15 March 2018.

\bibitem[{{Business Insider}(2017)}]{spacex}
{Business Insider}, \enquote{Elon Musk is about to launch the first of 11,925
  proposed SpaceX internet satellites — more than all spacecraft that orbit
  Earth today,}
  \url{http://www.businessinsider.com/spacex-starlink-microsat-launch-global-internet-2018-2},
  2017.
\newblock Accessed 16 February 2018.

\bibitem[{{Space News}(2017)}]{fccmega}
{Space News}, \enquote{FCC gets five new applications for non-geostationary
  satellite constellations,}
  \url{http://spacenews.com/fcc-gets-five-new-applications-for-non-geostationary-satellite-constellations/},
  2017.
\newblock Accessed 05 February 2018.

\bibitem[{Diekelman(1998)}]{diekelman1998design}
Diekelman, D., \enquote{Design guidelines for post-2000 constellations,}
  \emph{Mission Design \& Implementation of Satellite Constellations},
  Springer, 1998, pp. 11--21.
\newblock Doi: \url{https://doi.org/10.1007/978-94-011-5088-0_2}.

\bibitem[{{NASA Space Flight}(2017)}]{iridium}
{NASA Space Flight}, \enquote{Iridium marks new satellite network, 20 healthy
  satellites and 55 more to launch,}
  \url{https://www.nasaspaceflight.com/2017/07/iridium-satellite-network-55-more/},
  2017.
\newblock Accessed 29 January 2018.

\bibitem[{Palmerini(1998)}]{palmerini1998hybrid}
Palmerini, G.~B., \enquote{Hybrid Configurations for Satellite Constellations,}
  \emph{Mission Design \& Implementation of Satellite Constellations},
  Springer, 1998, pp. 81--89.
\newblock Doi: \url{https://doi.org/10.1007/978-94-011-5088-0_7}.

\bibitem[{{OneWeb}(2017)}]{onewebBuilt}
{OneWeb}, \enquote{OneWeb satellites breaks ground on the world’s first
  state-of-the-art high-volume satellite manufacturing facility,}
  \url{http://www.oneweb.world/press-releases/2017/oneweb-satellites-breaks-ground-on-the-worlds-first-state-of-the-art-high-volume-satellite-manufacturing-facility/},
  2017.
\newblock Accessed 25 January 2018.

\bibitem[{{Space Flight Now}(2015)}]{onewebLaunch}
{Space Flight Now}, \enquote{OneWeb launch deal called largest commercial
  rocket buy in history,}
  \url{https://spaceflightnow.com/2015/07/01/oneweb-launch-deal-called-largest-commercial-rocket-buy-in-history/},
  2015.
\newblock Accessed 29 January 2018.

\bibitem[{Lang and Adams(1998)}]{lang1998comparison}
Lang, T.~J., and Adams, W.~S., \enquote{A comparison of satellite
  constellations for continuous global coverage,} \emph{Mission Design \&
  Implementation of Satellite Constellations}, Springer, 1998, pp. 51--62.
\newblock Doi: \url{https://doi.org/10.1007/978-94-011-5088-0_5}.

\bibitem[{Lansard and Palmade(1998)}]{lansard1998satellite}
Lansard, E., and Palmade, J.-L., \enquote{Satellite constellation design:
  Searching for global cost-efficiency trade-offs,} \emph{Mission Design \&
  Implementation of Satellite Constellations}, Springer, 1998, pp. 23--31.
\newblock Doi: \url{https://doi.org/10.1007/978-94-011-5088-0_3}.

\bibitem[{Palmade et~al.(1998)Palmade, Frayssinhes, Martinot, and
  Lansard}]{palmade1998skybridge}
Palmade, J.-L., Frayssinhes, E., Martinot, V., and Lansard, E., \enquote{The
  SkyBridge constellation design,} \emph{Mission Design \& Implementation of
  Satellite Constellations}, Springer, 1998, pp. 133--140.
\newblock Doi: \url{https://doi.org/10.1007/978-94-011-5088-0_12}.

\bibitem[{Cornara et~al.(1999)Cornara, Beech, Bell{\'o}-Mora, and Martinez~de
  Aragon}]{cornara1999satellite}
Cornara, S., Beech, T., Bell{\'o}-Mora, M., and Martinez~de Aragon, A.,
  \enquote{Satellite constellation launch, deployment, replacement and
  end-of-life strategies,} \emph{13th Annual AIAA/USU Conference on Small
  Satellites}, 1999.

\bibitem[{Ereau and Saleman(1996)}]{ereau1996modeling}
Ereau, J.-F., and Saleman, M., \enquote{Modeling and simulation of a satellite
  constellation based on Petri nets,} \emph{Reliability and Maintainability
  Symposium, 1996 Proceedings. International Symposium on Product Quality and
  Integrity., Annual}, IEEE, 1996, pp. 66--72.
\newblock Doi: \url{https://doi.org/10.1109/RAMS.1996.500644}.

\bibitem[{Sumter(2003)}]{sumter2003optimal}
Sumter, B.~R., \enquote{Optimal replacement policies for satellite
  constellations,} Tech. rep., Air Force Inst. Technol., Wright-Patterson Air
  Force Base, 2003.

\bibitem[{Kelley and Dessouky(2004)}]{kelley2004minimizing}
Kelley, C., and Dessouky, M., \enquote{Minimizing the cost of availability of
  coverage from a constellation of satellites: evaluation of optimization
  methods,} \emph{Systems engineering}, Vol.~7, No.~2, 2004, pp. 113--122.
\newblock Doi: \url{https://doi.org/10.1002/sys.10059}.

\bibitem[{Dishon and Weiss(1966)}]{dishon1966communications}
Dishon, M., and Weiss, G.~H., \enquote{A communications satellite replenishment
  policy,} \emph{Technometrics}, Vol.~8, No.~3, 1966, pp. 399--409.
\newblock Doi: \url{https://doi.org/10.1080/00401706.1966.10490373}.

\bibitem[{G{\"u}m{\"u}s and G{\"u}neri(2007)}]{gumus2007multi}
G{\"u}m{\"u}s, A.~T., and G{\"u}neri, A.~F., \enquote{Multi-echelon inventory
  management in supply chains with uncertain demand and lead times: literature
  review from an operational research perspective,} \emph{Proceedings of the
  Institution of Mechanical Engineers, Part B: Journal of Engineering
  Manufacture}, Vol. 221, No.~10, 2007, pp. 1553--1570.
\newblock Doi: \url{https://doi.org/10.1243/09544054JEM889}.

\bibitem[{Ganeshan(1999)}]{ganeshan1999managing}
Ganeshan, R., \enquote{Managing supply chain inventories: A multiple retailer,
  one warehouse, multiple supplier model,} \emph{International Journal of
  Production Economics}, Vol.~59, No.~1, 1999, pp. 341--354.
\newblock Doi: \url{https://doi.org/10.1016/S0925-5273(98)00115-7}.

\bibitem[{Deuermeyer and Schwarz(1981)}]{schwarz1981}
Deuermeyer, B., and Schwarz, L.~B., \enquote{A model for the analysis of system
  service level in warehouse/retailer distribution systems: the identical
  retailer case,} \emph{Multi-level production/inventory control systems:
  theory and practice}, North-Holland, 1981.

\bibitem[{Costantino et~al.(2013)Costantino, Di~Gravio, and
  Tronci}]{costantino2013multi}
Costantino, F., Di~Gravio, G., and Tronci, M., \enquote{Multi-echelon,
  multi-indenture spare parts inventory control subject to system availability
  and budget constraints,} \emph{Reliability Engineering \& System Safety},
  Vol. 119, 2013, pp. 95--101.
\newblock Doi: \url{https://doi.org/10.1016/j.ress.2013.05.006}.

\bibitem[{Caglar et~al.(2004)Caglar, Li, and Simchi-Levi}]{caglar2004two}
Caglar, D., Li, C.-L., and Simchi-Levi, D., \enquote{Two-echelon spare parts
  inventory system subject to a service constraint,} \emph{IIE Transactions},
  Vol.~36, No.~7, 2004, pp. 655--666.
\newblock Doi: \url{https://doi.org/10.1080/07408170490278265}.

\bibitem[{Prussing and Conway(1993)}]{prussing1993orbital}
Prussing, J., and Conway, B., \emph{Orbital Mechanics}, Oxford University
  Press, 1993.

\bibitem[{{Walker}(1984)}]{1984JBIS...37..559W}
{Walker}, J.~G., \enquote{{Satellite constellations},} \emph{Journal of the
  British Interplanetary Society}, 1984, pp. 559--572.

\bibitem[{Jensen and Bard(2003)}]{jensen2003operations}
Jensen, P., and Bard, J., \emph{Operations Research Models and Methods}, Wiley,
  2003.

\bibitem[{Collopy(2003)}]{collopy2003assigning}
Collopy, P., \enquote{Assigning value to reliability in satellite
  constellations,} \emph{AIAA Space 2003 Conference \& Exposition}, 2003, p.
  6214.
\newblock Doi: \url{https://doi.org/10.2514/6.2003-6214}.

\bibitem[{{\c{C}}inlar(1975)}]{ccinlar1975exceptional}
{\c{C}}inlar, E., \enquote{Markov renewal theory: A survey,} \emph{Management
  Science}, Vol.~21, No.~7, 1975, pp. 727--752.
\newblock Doi: \url{https://doi.org/10.1287/mnsc.21.7.727}.

\bibitem[{Schwarz et~al.(1985)Schwarz, Deuermeyer, and
  Badinelli}]{schwarz1985fill}
Schwarz, L.~B., Deuermeyer, B.~L., and Badinelli, R.~D., \enquote{Fill-rate
  optimization in a one-warehouse N-identical retailer distribution system,}
  \emph{Management Science}, Vol.~31, No.~4, 1985, pp. 488--498.
\newblock Doi: \url{https://doi.org/10.1287/mnsc.31.4.488}.

\bibitem[{{Virgin Orbit}({})}]{virgin}
{Virgin Orbit}, \url{https://virginorbit.com/}, {}.
\newblock Accessed 20 February 2018.

\bibitem[{Stein(1987)}]{stein1987large}
Stein, M., \enquote{Large sample properties of simulations using Latin
  hypercube sampling,} \emph{Technometrics}, Vol.~29, No.~2, 1987, pp.
  143--151.
\newblock Doi: \url{https://doi.org/10.1080/00401706.1987.10488205}.

\bibitem[{Dubos et~al.(2010)Dubos, Castet, and Saleh}]{dubos2010statistical}
Dubos, G.~F., Castet, J.-F., and Saleh, J.~H., \enquote{Statistical reliability
  analysis of satellites by mass category: Does spacecraft size matter?}
  \emph{Acta Astronautica}, Vol.~67, No.~5, 2010, pp. 584--595.
\newblock Doi: \url{https://doi.org/10.1016/j.actaastro.2010.04.017}.

\bibitem[{Erlank and Bridges(2015)}]{erlank2015multicellular}
Erlank, A.~O., and Bridges, C.~P., \enquote{A multicellular architecture
  towards low-cost satellite reliability,} \emph{Adaptive Hardware and Systems
  (AHS), 2015 NASA/ESA Conference on}, IEEE, 2015, pp. 1--8.
\newblock Doi: \url{https://doi.org/10.1109/AHS.2015.7231152}.

\bibitem[{Castet and Saleh(2009)}]{castet2009satellite}
Castet, J.-F., and Saleh, J.~H., \enquote{Satellite and satellite subsystems
  reliability: Statistical data analysis and modeling,} \emph{Reliability
  Engineering \& System Safety}, Vol.~94, No.~11, 2009, pp. 1718--1728.
\newblock Doi: \url{https://doi.org/10.1016/j.ress.2009.05.004}.

\bibitem[{Castet and Saleh(2010)}]{castet2010beyond}
Castet, J.-F., and Saleh, J.~H., \enquote{Beyond reliability, multi-state
  failure analysis of satellite subsystems: A statistical approach,}
  \emph{Reliability Engineering \& System Safety}, Vol.~95, No.~4, 2010, pp.
  311--322.
\newblock Doi: \url{https://doi.org/10.1016/j.ress.2009.11.001}.

\bibitem[{Soyuz-2({})}]{soyouz}
Soyuz-2, \url{https://en.wikipedia.org/wiki/Soyuz-2#Launch_history}, {}.
\newblock Accessed 12 June 2017.

\bibitem[{SpaceX({})}]{falcon}
SpaceX, \enquote{{Launch Manifest},} \url{http://www.spacex.com/missions}, {}.
\newblock Accessed 06 March 2018.

\end{thebibliography}

\end{document}